# High Dimensional Inference in Partially Linear Models

Ying Zhu,* Zhuqing Yu,† and Guang Cheng‡

**Abstract**

We propose two semiparametric versions of the debiased Lasso procedure for the model $Y_i = X_i\beta_0 + g_0(Z_i) + \varepsilon_i$, where $\beta_0$ is high dimensional but sparse (exactly or approximately). Both versions are shown to have the same asymptotic normal distribution and do not require the minimal signal condition for statistical inference of any component in $\beta_0$. Our method also works when $Z_i$ is high dimensional provided that the function classes $\mathbb{E}(X_{ij}|Z_i)$s and $\mathbb{E}(Y_i|Z_i)$ belong to exhibit certain sparsity features, e.g., a sparse additive decomposition structure. We further develop a simultaneous hypothesis testing procedure based on multiplier bootstrap. Our testing method automatically takes into account of the dependence structure within the debiased estimates, and allows the number of tested components to be exponentially high.

Running title: High Dimensional Inference in Partial Linear Model

## 1 Introduction

Semiparametric regression is a longstanding statistical tool that leverages the flexibility of nonparametric models while avoiding the "curse of dimensionality" (Bickel et al, 1998). A leading example of semiparametric regression models is the partially linear regression

$$Y_i = X_i\beta_0 + g_0(Z_i) + \varepsilon_i, \ i = 1, ..., n. \tag{1}$$

In (1), $\beta_0 \in \mathbb{R}^p$ is an unknown vector and $g_0$ is an unknown function; $X := (X_i)_{i=1}^n \in \mathbb{R}^{n \times p}$ and $Z := (Z_i)_{i=1}^n \in \mathbb{R}^{n \times d}$ are observed covariates ($X_i$ and $Z_i$ denote the $i$th row of $X$ and $Z$, respectively), $Y := (Y_i)_{i=1}^n \in \mathbb{R}^n$ are the response variables, $\varepsilon = (\varepsilon_i)_{i=1}^n \in \mathbb{R}^n$ is a noise vector with $\mathbb{E}(\varepsilon_i) = 0$ and $\mathbb{E}(\varepsilon_i^2) = 1$, and independent of $(X, Z)$. Throughout the paper, we assume that the data $\{X_i, Z_i, Y_i\}_{i=1}^n$ are $i.i.d.$. The goal of this paper is to establish statistical inference results, e.g., confidence intervals and hypothesis testing, for the high dimensional component $\beta_0$ in presence of the nuisance function $g_0$. In particular, we assume that $p \geq n$ and $\beta_0$ exhibits sufficient sparsity (meaning that the ordered coefficients in $\beta_0$ decay sufficiently fast). Our method also works when $Z_i$ is high dimensional ($d \geq n$) provided that the function classes $\mathbb{E}(X_{ij}|Z_i)$s and $\mathbb{E}(Y_i|Z_i)$ belong to exhibit certain sparsity features, e.g., a sparse additive decomposition structure as defined in Raskutti, et al. (2012).

For statistical inference of $\beta_0$ in (1), existing results mainly focus on the regime where $p$ increases with $n$ but smaller than $n$, for example, Li and Liang (2008), Xie and Huang (2009), and Cheng, et al. (2015). Sherwood and Wang (2016) allow $p \geq n$ but require the minimal signal condition. Such results therefore suffer the problems arising from the nonuniformity of limit theory. Recently, Javanmard and Montanari (2014), van de Geer, et al. (2014), and Zhang and Zhang (2014) have proposed the debiased Lasso for high dimensional linear models. These estimators are non-sparse, have a limiting normal distribution, and do not require the minimal signal condition. For the linear model $Y = X\beta_0 + \varepsilon$, given an initial Lasso estimate $\hat{\beta}$ of $\beta_0$, the debiased Lasso adds a correction term to $\hat{\beta}_j$ (the $j$th component of $\hat{\beta}$) to remove the bias introduced by regularization. In particular, the correction term takes the form of

$$\hat{\Gamma}_j \frac{1}{n} X^T \left( Y - X\hat{\beta} \right). \tag{2}$$

*Research Associate. Department of Economics. Michigan State University. zip: 48824. yzhu@msu.edu
†Graduate Student. Department of Statistics. Purdue University.
‡Professor of Statistics. Department of Statistics. Purdue University. Research Sponsored by NSF CAREER Award DMS-1151692, DMS-1418042, and Office of Naval Research.

In (2), $n^{-1} X^T \left( Y - X\hat{\beta} \right)$ is the sample analogue of the population score function $\mathbb{E}\left( X_i^T (Y_i - X_i \beta_0) \right)$; $\hat{\Gamma}_j$ denotes the $j$th row of $\hat{\Gamma}$ where $\hat{\Gamma}$ is an approximate inverse of $n^{-1} X^T X$, whose population counterpart is $\mathbb{E}\left( X_i^T X_i \right)$. In our model (1), additional bias arises due to the presence of $g_0$; consequently, the standard debiased Lasso cannot rid of the effect from $g_0$ and thus will not have a limiting distribution centered around zero. Instead, we propose two modified versions of the debiased Lasso estimators for $\beta_0$. Both versions are shown to be asymptotically unbiased for $\beta_0$, have the same limiting (normal) distribution, and do not require the minimal signal condition.

Our modified debiased Lasso estimators use a "nonparametric projection" strategy to remove the impact of $g_0$ in (1). Such a strategy has been used in the semiparametric inference literature where $p$ is assumed to be small relative to $n$ (e.g., Robinson, 1988; Donald and Newey, 1994; Liang and Li, 2009). To be more specific, by taking the conditional expectations of the left side and the right side of (1) with respect to $Z_i$, we obtain

$$\mathbb{E}\left( Y_i | Z_i \right) = \mathbb{E}\left( X_i | Z_i \right) \beta_0 + g_0(Z_i) \tag{3}$$

where we exploit the fact that $\mathbb{E}\left( \varepsilon_i | Z_i \right) = 0$. Subtracting $\mathbb{E}\left( Y_i | Z_i \right)$ from $Y_i$ and $\mathbb{E}\left( X_i | Z_i \right) + g_0(Z_i)$ from $X_i + g_0(Z_i)$ in (1) yields

$$\tilde{Y}_i = \tilde{X}_i \beta_0 + \varepsilon_i \tag{4}$$

where $\tilde{Y}_i := Y_i - \mathbb{E}\left( Y_i | Z_i \right)$, $\tilde{X}_{ij} := X_{ij} - \mathbb{E}\left( X_{ij} | Z_i \right)$ and $\tilde{X}_i := \left( \tilde{X}_{ij} \right)_{j=1}^p$ (which is a $p$-dimensional row vector).

Relating (4) to the linear model $Y_i = X_i \beta_0 + \varepsilon_i$, given nonparametric surrogates $\hat{Y}_i := Y_i - \hat{\mathbb{E}}\left( Y_i | Z_i \right)$ of $\tilde{Y}_i$ and $\hat{X}_{ij} := X_{ij} - \hat{\mathbb{E}}\left( X_{ij} | Z_i \right)$ of $\tilde{X}_{ij}$ ($j = 1, ..., p$), we simply replace $Y_i$ with $\hat{Y}_i$, $X_i$ with the row vector $\hat{X}_i := \left( \hat{X}_{ij} \right)_{j=1}^p$, and $\hat{\Gamma}_j$ with the $j$th row ($\hat{\Theta}_j$) of an approximate inverse (denoted as $\hat{\Theta}$) of $n^{-1} \sum_{i=1}^n \hat{X}_i^T \hat{X}_i$ in (2). This yields our first semiparametric version of the debiased procedure

$$\hat{b}_j := \hat{\beta}_j + \hat{\Theta}_j \frac{1}{n} \hat{X}^T \left( \hat{Y} - \hat{X}\hat{\beta} \right), \tag{5}$$

where $\hat{\beta}$ is an initial estimate of $\beta_0$. Alternatively, by noting that $n^{-1} \hat{X}^T \left( \hat{Y} - \hat{X}\hat{\beta} \right)$ in (5) is simply the sample analogue of the population score function $\mathbb{E}\left( \tilde{X}_i^T \varepsilon_i \right)$, we arrive at our second debiased procedure

$$\tilde{b}_j := \hat{\beta}_j + \hat{\Theta}_j \frac{1}{n} \hat{X}^T \left( Y - X\hat{\beta} - \hat{g} \right), \tag{6}$$

where $\hat{g}$ is an estimate of $g_0$.

We provide theoretical implications on the impact of the estimation errors associated with the $p$ nonparametric surrogates $\hat{\mathbb{E}}\left( X_{ij} | Z_i \right)$s in our modified debiased Lasso procedures when each of $\hat{\mathbb{E}}(X_{ij}|Z_i)$s concerns a large family of (regularized) nonparametric least squares estimators. These implications also hold true for the surrogate $\hat{\mathbb{E}}(Y_i|Z_i)$ (which matters to (5)) and the surrogate $\hat{g}(Z_i)$ (which matters to (6)). After careful theoretical analysis, we find that if the error of the nonparametric regression *per se* (with respect to the prediction norm) is $O_p(r_n)$, it only contributes $O_p(r_n^2)$ in the asymptotic expansions of $\hat{b}_j - \beta_{0j}$ and $\tilde{b}_j - \beta_{0j}$ for any $j = 1, ..., p$, where $r_n$ is related to the optimal rate for the nonparametric regression. This result implies that even with $p$ much larger than $n$ (and/or with the dimension $d$ of $Z_i$ much larger than $n$), the limiting distribution of our modified debiased estimators for any individual component in $\beta_0$ may behave as if the unknown conditional expectations $\mathbb{E}(X_{ij}|Z_i)$s and $\mathbb{E}(Y_i|Z_i)$ as well as the unknown function $g_0(Z_i)$ were known.

This theoretical finding motivates us to consider a multiplier-bootstrap-based *simultaneous* hypothesis testing procedure for any sub-vector of $\beta_0$. This extends the method developed by Zhang and Cheng (2017) from linear regressions to more flexible partially linear regressions. Our simultaneous testing procedure automatically takes into account of the dependence structure within our debiased estimators, and allows the number of tested components to be exponentially high.

We illustrate the theoretical finding with four specific examples in terms of $\dim(Z_i)$ and the function class that $\mathbb{E}(X_{ij}|Z_i)$s and $\mathbb{E}(Y_i|Z_i)$ belong to. With regard to the specific forms of $\hat{\mathbb{E}}(X_{ij}|Z_i)$s and $\hat{\mathbb{E}}(Y_i|Z_i)$, several modern techniques for the projection step are considered and the rates achieved by these practical procedures are compared with the theoretical results. The techniques discussed in the paper include the smoothing splines estimator in Sobolev balls, the Lasso (Tibshirani, 1996) and Slope (Su and Candès, 2016) in sparse linear regression models, the (square-root) Lasso (Belloni, et al., 2014) in rearranged Sobolev balls, and the $l_1$−regularized kernel ridge regression (Raskutti et al., 2012) in sparse additive models.

When constructing an approximate inverse $\hat{\Theta}$ of $n^{-1} \sum_{i=1}^{n} \hat{X}_i^T \hat{X}_i$ in (5) and (6), we adopt the nodewise regression method proposed by van de Geer, et al. (2014). Since our analysis involves establishing $\left\| \hat{\Theta}_j - \Theta_j \right\|_1 = o_p(1)$, as in van de Geer et al. (2014), we require a sparsity condition on the inverse $\Theta = \Sigma^{-1}$ of the population Hessian $\Sigma := \mathbb{E}\left( \tilde{X}_i^T \tilde{X}_i \right)$. This fact renders the method in Javanmard and Montanari (2014) for constructing $\hat{\Theta}$ inapplicable as their approach is only valid for fixed $X$, while our analysis accounts for the randomness in $X$ and the estimation errors in $\hat{\mathbb{E}}(X_{ij}|Z_i)$s. Furthermore, our analysis relaxes the exact sparsity of $\Theta_j$ (assumed in most literature including van de Geer, et al., 2014) to accommodate for approximate sparsity which permits all the entries in $\Theta_j$ to be non-zero as long as they decay sufficiently fast. This extension provides a more realistic interpretation of most practical problems.

The rest of the paper is organized as follows. Section 2 presents the detailed construction of the modified debiased estimators for $\beta_0$ and the simultaneous testing procedure. Section 3 establishes the main theoretical results. All technical details are deferred to Section 4. Section 5 evaluates the performance of our methods with simulation experiments.

**Notation**. The $l_q$−norm of a $p$−dimensional vector $\Delta$ is denoted by $\|\Delta\|_q$ for $1 \leq q \leq \infty$. For a matrix $H \in \mathbb{R}^{p_1 \times p_2}$, write $\|H\|_\infty := \max_{i,j} |H_{ij}|$ to be the elementwise $l_\infty$−norm of $H$. Let $H_j$ denote the $j$th row of $H$. The $\mathcal{L}^2(\mathbb{P}_n)$−norm of the vector $\Delta := \{\Delta(X_i)\}_{i=1}^n$, denoted by $\|\Delta\|_n$, is given by $\left[ \frac{1}{n} \sum_{i=1}^n (\Delta(X_i))^2 \right]^{\frac{1}{2}}$. For functions $f(n)$ and $g(n)$, write $f(n) \gtrsim g(n)$ to mean that $f(n) \geq cg(n)$ for a universal constant $c \in (0, \infty)$ and similarly, $f(n) \lesssim g(n)$ to mean that $f(n) \leq c' g(n)$ for a universal constant $c' \in (0, \infty)$, and $f(n) \asymp g(n)$ when $f(n) \gtrsim g(n)$ and $f(n) \lesssim g(n)$ hold simultaneously. Also denote $\max\{a, b\}$ by $a \vee b$ and $\min\{a, b\}$ by $a \wedge b$. As a general rule for this paper, all the $c \in (0, \infty)$ constants denote positive universal constants. The specific values of these constants may change from place to place.

## 2 Main methodology

In this section we discuss the construction of $\hat{b}_j$ and $\tilde{b}_j$ in detail. Note that both (5) and (6) require estimators for the conditional expectations, an initial estimator $\hat{\beta}$ for $\beta_0$ (and an estimator $\hat{g}$ for $g_0$ in $\tilde{b}_j$), and also an approximate inverse $\hat{\Theta}$ for $n^{-1}\sum_{i=1}^n \hat{X}_i^T \hat{X}_i$. We first discuss how to obtain these aforementioned quantities. Given $\hat{b}_j$s and $\tilde{b}_j$s, we then present the simultaneous inference procedure.

**Estimators for the conditional expectations**

For either (5) or (6), we need to estimate the conditional expectations $\mathbb{E}(X_{ij}|Z_i)$s ($j = 1,...,p$). This step is easily parallelable as it involves solving $p$ independent subproblems and each subproblem can be in general solved with an efficient algorithm. In contrast with (6), (5) does not require an estimate of $g_0$ but an estimate of $\mathbb{E}(Y_i|Z_i)$. Estimating conditional expectations is widely studied in the literature on nonparametric methods. For the purpose of this paper, global properties of the nonparametric estimators $\hat{\mathbb{E}}(X_{ij}|Z_i)$s and $\hat{\mathbb{E}}(Y_i|Z_i)$ are the key to our analysis of the debiased procedures and therefore, we focus on the following least squares estimators

$$\hat{f}_j \in \arg\min_{f_j \in \mathcal{F}_j} \left\{ \frac{1}{2n} \sum_{i=1}^n (w_{ij} - f_j(z_i))^2 \right\}, \tag{7}$$

where $w_{i0} = y_i$ and $w_{ij} = x_{ij}$ for $j = 1,...,p$. Denote $\hat{f}_0(Z_i)$ as $\hat{\mathbb{E}}(Y_i|Z_i)$ and $\hat{f}_j(Z_i)$ as $\hat{\mathbb{E}}(X_{ij}|Z_i)$.

A nonparametric regression problem like (7) is a standard setup in many modern statistics books (e.g., van de Geer, 2000; Wainwright, 2015). Examples of (7) include linear regression, sparse linear regression, series projection, convex regression, Lipschitz Isotonic regression, and kernel ridge regression (KRR). In the case of KRR, we restrict $\mathcal{F}_j$ in (7) to be a compact subset of a reproducing kernel Hilbert space (RKHS) $\mathcal{H}$, equipped with a norm $\|\cdot\|_\mathcal{H}$; (7)[1] can then be reformulated in its Lagrangian form

$$\hat{f}_j \in \arg\min_{f_j \in \mathcal{H}} \left\{ \frac{1}{2n} \sum_{i=1}^n (w_{ij} - f_j(z_i))^2 + \mu_j \|f_j\|_\mathcal{H}^2 \right\} \tag{8}$$

where $\mu_j > 0$ is a regularization parameter. In particular, smoothing spline estimators can be viewed as special cases of KRR.

**Initial estimators for $\beta_0$ and $g_0$**

In a semiparametric regression model like (1), Zhu (2017) covers a wide spectrum of function classes that the nonparametric component $g_0(\cdot)$ may belong to and provides a general nonasymptotic theory for estimating $\beta_0$ and $g_0$. The estimators $\hat{\beta}$ and $\hat{g}$ in Zhu (2017) can be used as initial estimators in (5)-(6). Given the way $\hat{\beta}$ is obtained in Zhu (2017), the estimated conditional expectations, $\hat{f}_j$s, come in handy as byproducts (therefore, separate estimations for the conditional expectations are not needed in the construction of $\hat{b}_j$ or $\tilde{b}_j$).

For special cases where $Z_i$ has a low dimension and $g_0$ belongs to the $m$th order Sobolev ball $\mathcal{S}^m$, other estimators for $\beta_0$ and $g_0$ are also available (Müller and van de Geer, 2015; Yu, et al. 2017) and take the following form

$$(\hat{\beta}, \hat{g}) \in \arg\min_{\beta \in \mathbb{R}^p, g \in \mathcal{S}^m} \left\{ \frac{1}{2n} \sum_{i=1}^n (Y_i - X_i\beta - g(Z_i))^2 + \mu J^2(g) + \lambda \|\beta\|_1 \right\}, \tag{9}$$

---
[1] To be more specific, we let $\mathcal{F}_j$ be a ball of radius $R$ in the norm $\|\cdot\|_\mathcal{H}$ and assume $R \leq 1$ throughout the asymptotic analysis to avoid carrying "$R$"s around.

where $\mu, \lambda > 0$ are regularization parameters, $J^2(g) := \int_0^1 (g^m(z))^2 dz$ and $Z \sim \text{Uniform}(0, 1)$.

Due to the intractable limiting distribution of Lasso type estimators, these aforementioned papers do not provide any distributional results for their proposed estimators. In Section 3, we take the debiased versions, (5) and (6), of these aforementioned initial estimators and establish the asymptotic normality of individual components in the debiased estimators.

**Estimator for the inverse of the population Hessian**

Given the estimates $\hat{Y}_i$ of $\tilde{Y}_i$ and $\hat{X}_i$ of $\tilde{X}_i$ via (7), we obtain an approximate inverse $\hat{\Theta}$ of $n^{-1} \sum_{i=1}^n \hat{X}_i^T \hat{X}_i$ using the nodewise regression method proposed by van de Geer, et al. (2014). Since our analysis involves establishing $\left\|\hat{\Theta}_j - \Theta_j\right\|_1 = o_p(1)$, as in van de Geer, et al. (2014), we require a sparsity condition on the inverse $\Theta = \Sigma^{-1}$ of the population Hessian $\Sigma := \mathbb{E}\left(\tilde{X}_i^T \tilde{X}_i\right)$. Lack of sparsity in the off-diagonal elements of $\Theta$ will cause remainder terms like $\left(\hat{\Theta}_j - \Theta_j\right) \frac{1}{\sqrt{n}} \tilde{X}^T \varepsilon$ in the asymptotic expansions[2] of $\sqrt{n}\left(\hat{b}_j - \beta_{0j}\right)$ or $\sqrt{n}\left(\tilde{b}_j - \beta_{0j}\right)$ to diverge and the resulting limiting distribution may not be well-behaved for any practical purpose. This fact renders the method in Javanmard and Montanari (2014) for constructing $\hat{\Theta}$ inapplicable as their approach is only valid for fixed $X$, while our analysis accounts for the randomness in $X$ and the estimation errors in $\hat{\mathbb{E}}(X_{ij}|Z_i)$s.

To apply the nodewise regression method in our context, for each $1 \leq j \leq p$, let us define

$$\hat{\pi}_j \in \arg\min_{\tilde{\pi}_j \in \mathbb{R}^{p-1}} \left\{ \frac{1}{n} \left\|\hat{X}_j - \hat{X}_{-j} \tilde{\pi}_j\right\|_2^2 + \lambda_j \|\tilde{\pi}_j\|_1 \right\}, \tag{10}$$

where $\hat{X}_{-j}$ denotes the submatrix of $\hat{X}$ without the $j$th column. Let

$$\hat{M} := \begin{pmatrix} 1 & -\hat{\pi}_{1,2} & \cdots & -\hat{\pi}_{1,p} \\ -\hat{\pi}_{2,1} & 1 & \cdots & -\hat{\pi}_{2,p} \\ \vdots & \vdots & \ddots & \vdots \\ -\hat{\pi}_{p,1} & -\hat{\pi}_{p,2} & \cdots & 1 \end{pmatrix}.$$

Based on $\hat{\pi}_j := \left\{\hat{\pi}_{j,j'}; j' \neq j\right\}$, for each $1 \leq j \leq p$, we compute

$$\hat{\tau}_j^2 := \frac{1}{n} \left\|\hat{X}_j - \hat{X}_{-j} \hat{\pi}_j\right\|_2^2 + \lambda_j \|\hat{\pi}_j\|_1$$

and write

$$\hat{T}^2 := \text{diag}\left(\hat{\tau}_1^2, ..., \hat{\tau}_p^2\right).$$

Finally, we define $\hat{\Theta} := \hat{T}^{-2} \hat{M}$.

For later presentations of the theoretical results, we also introduce the population counterparts of the above quantities: let $\pi_j$ be the population regression coefficients of $\tilde{X}_{ij}$ on $\tilde{X}_{i,-j} = \left\{\tilde{X}_{ij'}; j' \neq j\right\}$ and

$$M := \begin{pmatrix} 1 & -\pi_{1,2} & \cdots & -\pi_{1,p} \\ -\pi_{2,1} & 1 & \cdots & -\pi_{2,p} \\ \vdots & \vdots & \ddots & \vdots \\ -\pi_{p,1} & -\pi_{p,2} & \cdots & 1 \end{pmatrix},$$

$$T^2 := \text{diag}\left(\tau_1^2, ..., \tau_p^2\right),$$

---

[2]See (21)-(22) for more detail.

such that $\tau_j^2 := \mathbb{E}\left[\left(\tilde{X}_{ij} - \tilde{X}_{i,-j}\pi_j\right)^2\right]$ for $j = 1, ..., p$.

**Simultaneous inference**

From a practical viewpoint, conducting simultaneous inference for a collection of parameters in high-dimensional models may be of greater interest to researchers than inference of a single parameter. To be more specific, suppose we are interested in testing the hypothesis:

$$H_{0,G} : \beta_{0j} = \tilde{\beta}_j \ \forall j \in G \subseteq \{1, 2, ..., p\}$$

against the alternative $H_{a,G} : \beta_{0j} \neq \tilde{\beta}_j$ for some $j \in G$. In particular, we allow $|G| \geq n$. Zhang and Cheng (2017) develop a bootstrap-assisted procedure to conduct simultaneous inference in sparse linear models. Here we propose similar test statistics

$$\begin{aligned}
T_G &= \max_{j \in G} \sqrt{n}\left(\hat{b}_j - \tilde{\beta}_j\right), \\
\text{or, } T_G &= \max_{j \in G} \sqrt{n}\left(\tilde{b}_j - \tilde{\beta}_j\right),
\end{aligned}$$

and a multiplier bootstrap version

$$W_G = \max_{j \in G} \frac{1}{\sqrt{n}} \sum_{i=1}^n \hat{\Theta}_j \hat{X}_i^T \epsilon_i,$$

where $(\epsilon_i)_{i=1}^n$ are i.i.d. $\mathcal{N}(0, 1)$ random variables. The bootstrap critical value is then given by

$$c_G(\alpha) = \inf\{t \in \mathbb{R} : \mathbb{P}(W_G \leq t|Y, X, Z) \geq 1 - \alpha\}$$

for any user-defined $\alpha \in (0, 1)$. In the case where the variance $\sigma_\varepsilon^2$ of $\varepsilon_i$ in (1) is unknown, we can use

$$W_G = \max_{j \in G} \frac{1}{\sqrt{n}} \sum_{i=1}^n \hat{\Theta}_j \hat{X}_i^T \hat{\sigma}_\varepsilon \epsilon_i,$$

where $\hat{\sigma}_\varepsilon^2 = \frac{\sum_{i=1}^n (\hat{Y}_i - X_i\hat{\beta})^2}{n - \|\hat{\beta}\|_1}$ or $\hat{\sigma}_\varepsilon^2 = \frac{\sum_{i=1}^n (Y_i - X_i\hat{\beta} - \hat{g}(Z_i))^2}{n - \|\hat{\beta}\|_1}$ is an estimator for $\sigma_\varepsilon^2$.

## 3 Theoretical results

To make the key point of this paper, we first present the results for the case where $\beta_0$ and $\pi_j$ are assumed to be exactly sparse (Theorem 1). Our second result (Theorem 2 in Section 3.4) relaxes the exact sparsity assumptions and allows $\beta_0$ and $\pi_j$ to be approximately sparse. Note that the (exact or approximate) sparsity of $\pi_j$ along with the condition $\frac{1}{\tau_j^2} \precsim 1$ implies the sparsity of $\Theta_j = \left(\left[\mathbb{E}\left(\tilde{X}_i^T \tilde{X}_i\right)\right]^{-1}\right)_j$.

We begin with the following definitions. Let

$$\begin{aligned}
s_0 &:= |\{j : \beta_{0j} \neq 0\}|, \\
s_j &:= \left|\left\{j' \neq j, 1 \leq j' \leq p : \pi_{j,j'} \neq 0\right\}\right|.
\end{aligned}$$

To simplify our notations, we assume that $\|\beta_0\|_1 \precsim s_0$ and $\|\pi_j\|_1 \precsim s_j$ in Theorem 1 and Corollary 1. Recall $\mathcal{F}_j$ in (7); for notation simplicity, we assume $\mathcal{F}_j = \mathcal{F}$ from now on. Note that this restriction can be easily relaxed in our analysis. For any radius $\tilde{r}_n > 0$, we define the conditional *local complexity*

$$\mathcal{G}_n(\tilde{r}_n; \mathcal{F}) := \mathbb{E}_\xi \left[ \sup_{f \in \Omega(\tilde{r}_n; \mathcal{F})} \left| \frac{1}{n} \sum_{i=1}^n \xi_i f(Z_i) \right| \mid \{Z_i\}_{i=1}^n \right], \tag{11}$$

where $\xi_i$s are i.i.d. zero-mean sub-Gaussian variables with parameter at most $\sigma^\dagger$ and $\mathbb{E}(\xi_i | Z_i) = 0$ for all $i = 1, ..., n$, and

$$\Omega(\tilde{r}_n; \mathcal{F}) := \left\{ f \in \bar{\mathcal{F}} : \|f\|_n \leq \tilde{r}_n \right\},$$
$$\bar{\mathcal{F}} := \left\{ f = f' - f'' : f', f'' \in \mathcal{F} \right\}.$$

For any star-shaped class $\bar{\mathcal{F}}$ (that is, for any $f \in \bar{\mathcal{F}}$ and $\alpha \in [0, 1]$, $\alpha f \in \bar{\mathcal{F}}$), Lemma A8 in Section 4 guarantees that the function $t \mapsto \frac{\mathcal{G}_n(t; \mathcal{F})}{t}$ is non-increasing on the interval $(0, \infty)$. Therefore, there exists some large enough $\tilde{r}_n > 0$ that satisfies the *critical inequality*

$$\mathcal{G}_n(\tilde{r}_n; \mathcal{F}) \leq \frac{\tilde{r}_n^2}{2}; \tag{12}$$

moreover, (12) has a smallest positive solution $r_n$ (which we will refer to as the *critical radius*). In practice, determining the exact value of this *critical radius* can be difficult; fortunately, reasonable upper bounds on $r_n$ are often available. Here we describe two common methods from existing literature.

By a discretization argument and the Dudley's entropy integral, we may bound (11) by

$$c_0 \left( \frac{\sigma^\dagger}{\sqrt{n}} \int_0^{\tilde{r}_n} \sqrt{\log N_n(t; \Omega(\tilde{r}_n; \mathcal{F}))} dt + \tilde{r}_n^2 \right)$$

for some universal constant $c_0 > 0$, where $N_n(t; \Omega(\tilde{r}_n; \mathcal{F}))$ is the $t-$covering number of the set $\Omega(\tilde{r}_n; \mathcal{F})$. Let $\tilde{r}_n$ be a solution for

$$\frac{\sigma^\dagger}{\sqrt{n}} \int_0^{\tilde{r}_n} \sqrt{\log N_n(t; \Omega(\tilde{r}_n; \mathcal{F}))} dt \precsim \tilde{r}_n^2. \tag{13}$$

The resulting $\tilde{r}_n$ is known to yield an upper bound on the *critical radius* $r_n$ for (12) (see Lemma A9 for a formal statement); moreover, such bounds achieve sharp scaling on $r_n$ for a wide variety of function classes (see e.g., Barlett and Mendelson, 2002; Koltchinski, 2006; Wainwright, 2015).

When $\mathcal{F}$ is a ball of radius $R$ in the RKHS norm $\|\cdot\|_\mathcal{H}$, we let

$$\Omega(\tilde{r}_n; \mathcal{F}) := \left\{ f \in \bar{\mathcal{F}} : \|f\|_n \leq \tilde{r}_n, \|f\|_\mathcal{H} \leq 1 \right\}.$$

In this case, we can determine a good upper bound for $r_n$ using the result in Mendelson (2002) who shows that

$$\mathcal{G}_n(\tilde{r}_n; \mathcal{F}) \precsim \sigma^\dagger \sqrt{\frac{1}{n} \sum_{i=1}^n \min\{\tilde{r}_n^2, \tilde{\mu}_i\}}$$

where $\tilde{\mu}_1 \geq \tilde{\mu}_2 \geq ... \geq \tilde{\mu}_n \geq 0$ are the eigenvalues of the underlying kernel matrix for the KRR estimate. Consequently, we can solve for $\tilde{r}_{nj}$ via

$$\sigma^\dagger \sqrt{\frac{1}{n} \sum_{i=1}^n \min\{\tilde{r}_n^2, \tilde{\mu}_i\}} \precsim \tilde{r}_n^2.$$

This method above is known to yield $\tilde{r}_{nj}$ with sharp scaling for various choices of kernels.

## 3.1 The case of exact sparsity

We are now ready to establish our first main result (Theorem 1), which requires the following assumptions (in addition to those stated at the beginning of Section 1).

**Assumption 1.** For $j = 1, ..., p$, $\left(X_{ij}, \tilde{X}_{ij}, \tilde{Y}_i, \varepsilon_i\right)_{i=1}^n$ are i.i.d. sub-Gaussian variables with parameters at most $O(1)$.

**Assumption 2.** (i) For $\Sigma = \mathbb{E}\left(\tilde{X}_i^T \tilde{X}_i\right)$, the smallest eigenvalue $\Lambda_{\min}^2(\Sigma)$ is bounded away from 0, i.e., $\frac{1}{\Lambda_{\min}^2(\Sigma)} \precsim 1$; moreover, $\Lambda_{\max}^2(\Sigma) \precsim 1$. (ii) For $\Sigma_1 = \mathbb{E}\left[\mathbb{E}(X_i|Z_i)^T \mathbb{E}(X_i|Z_i)\right]$ (where $\mathbb{E}(X_i|Z_i) := (\mathbb{E}(X_{i1}|Z_i), ..., \mathbb{E}(X_{ip}|Z_i))$ is a $p$-dimensional row vector), $\Lambda_{\max}^2(\Sigma_1) = O(1)$.

**Remark.** Part (ii) of Assumption 2 is only used in the analysis for the second debiased estimator $\tilde{b}$ in (6).

**Assumption 3.** The conditional expectations $\mathbb{E}(Y_i|Z_i)$ and $\mathbb{E}(X_{ij}|Z_i)$ $(j = 1, ..., p)$ belong to $\mathcal{F}$. For any $f \in \bar{\mathcal{F}} = \left\{\tilde{f} = f' - f'' : f', f'' \in \mathcal{F}\right\}$ and $\alpha \in [0, 1]$, $\alpha f \in \bar{\mathcal{F}}$ (that is, $\bar{\mathcal{F}}$ is star-shaped).

**Remark.** This condition is needed to operationalize (5)-(6) which require estimators for the conditional expectations. In view of the following identity

$$g_0(Z_i) = \mathbb{E}(Y_i|Z_i) - \mathbb{E}(X_i|Z_i)\beta_0,$$

note that imposing conditions on the function class $\mathbb{E}(X_{ij}|Z_i)$s and $\mathbb{E}(Y_i|Z_i)$ belong to automatically restricts the function class $g_0$ belongs to.

**Assumption 4.** (i) The initial estimator $\hat{\beta}$ satisfies that $\left\|\hat{\beta} - \beta_0\right\|_2 = O_p\left(\sqrt{\frac{s_0 \log p}{n}}\right)$ and $\left\|\hat{\beta} - \beta_0\right\|_1 = O_p\left(s_0\sqrt{\frac{\log p}{n}}\right)$. (ii) $\left\|\mathbb{E}(X|Z)\left(\hat{\beta} - \beta_0\right)\right\|_n = O_p\left(\sqrt{\frac{s_0 \log p}{n}}\right)$ and the estimator $\hat{g}$ satisfies that $\|\hat{g} - g\|_n^2 = O_p\left((s_0 r_n^2) \vee \left(\frac{s_0 \log p}{n}\right)\right)$ where $r_n$ is the critical radius.

**Remark.** Part (ii) of Assumption 4 is only used in the analysis for the second debiased estimator $\tilde{b}$ in (6). Under mild conditions, the rates in Assumption 4 are satisfied by the initial estimators based on Zhu (2017). For special cases where $Z_i$ has a low dimension and $g_0$ belongs to the $m$th order Sobolev ball $\mathcal{S}^m$, the initial estimators based on Müller and van de Geer (2015) or Yu, et al. (2017) also satisfy the rate requirements in Assumption 4.

**Assumption 5.** $\left(\|\pi_j\|_1 \vee 1\right) r_n^2 = O\left(\sqrt{\frac{\log p}{n}}\right)$ where $r_n$ is the critical radius. Moreover,

$$\sqrt{n} s_j (s_0 \vee 1) \left(r_n^2 \vee \frac{\log p}{n}\right) = o(1),$$

$$\sqrt{n} \left\{\left[s_j^2 (s_0 \vee 1) \frac{\log p}{n}\right] \vee \left[s_j^3 (s_0 \vee 1) \left(r_n^2 \vee \frac{\log p}{n}\right) \sqrt{\frac{\log p}{n}}\right]\right\} = o(1),$$

$$\left[(s_j^2 \vee 1) \sqrt{\frac{\log p}{n}}\right] \vee \left[(s_j^3 \vee 1) \left(r_n^2 \vee \frac{\log p}{n}\right)\right] = o(1).$$

**Remark**. With some algebraic manipulations, we show that $\sqrt{n}\left(\hat{b}_j - \beta_{0j}\right) = \frac{1}{\sqrt{n}}\Theta_j \tilde{X}^T \varepsilon + REM$ and $\sqrt{n}\left(\tilde{b}_j - \beta_{0j}\right) = \frac{1}{\sqrt{n}}\Theta_j \tilde{X}^T \varepsilon + REM'$, where the leading term $\frac{1}{\sqrt{n}}\Theta_j \tilde{X}^T \varepsilon$ has an asymptotic normal distribution; see (21) and (22) for more details on the forms of $REM$ and $REM'$. Assumption 5 imposes requirements on the sparsity parameters $s_0$ and $s_j$s as well as the rates (reflected by $r_n$) for the auxiliary estimators $\hat{f}_j$s so that $REM = o(1)$ and $REM' = o(1)$.

**Assumption 6**. For all $l \neq j$, $j = 1, ..., p$, $\mathbb{E}\left[\frac{1}{n}\tilde{X}_l^T\left(\tilde{X}_j - \tilde{X}_{-j}\pi_j\right)\right] = O\left(\sqrt{\frac{\log p}{n}}\right)$ as $p \to \infty$ and $n \to \infty$.

**Remark**. For a general sub-Gaussian matrix, Assumption 6 is needed in order to derive the scaling for $\lambda_j$, whose choice in (10) depends on an upper bound for $\frac{1}{n}\left\|\tilde{X}_{-j}^T\left(\tilde{X}_j - \pi_j \tilde{X}_{-j}\right)\right\|_\infty$. Intuitively, Assumption 6 says that as the number of terms ($\tilde{X}_{i,-j}$) used to approximate $\tilde{X}_{ij}$ increases (that is, as $p \to \infty$), $\mathbb{E}\left[\frac{1}{n}\tilde{X}_l^T\left(\tilde{X}_j - \pi_j \tilde{X}_{-j}\right)\right] = o(1)$ provided $\frac{\log p}{n} = o(1)$. If $\mathbb{E}\left(\tilde{X}_{ij}|\tilde{X}_{i,-j}\right) = \tilde{X}_{i,-j}\pi_j$ (e.g., when $\tilde{X}_i$ is a normal vector), then $\mathbb{E}\left[\frac{1}{n}\tilde{X}_l^T\left(\tilde{X}_j - \pi_j \tilde{X}_{-j}\right)\right] = 0$ for all $l \neq j$. This special case is considered in Theorem 2.2 in van de Geer, et al. (2014).

**Theorem 1**. Under Assumptions 1-6, if we choose $\lambda_j \asymp \sqrt{\frac{\log p}{n}}$ uniformly in $j$ in (10), then

$$\frac{\sqrt{n}\left(\hat{b}_j - \beta_{0j}\right)}{\hat{\sigma}_j} \xrightarrow{D} \mathcal{N}(0, 1),$$

$$\frac{\sqrt{n}\left(\tilde{b}_j - \beta_{0j}\right)}{\hat{\sigma}_j} \xrightarrow{D} \mathcal{N}(0, 1),$$

where $\hat{\sigma}_j^2 = \hat{\Theta}_j \frac{\hat{X}^T\hat{X}}{n} \hat{\Theta}_j^T$, for each $j = 1, ..., p$.

Based on Theorem 1, Corollary 1 justifies the use of multiplier bootstrap in testing $H_{0,G}$ even when $|G|$ diverges.

**Corollary 1**. Suppose Assumptions 1-4, 6 hold and Assumption 5 is satisfied with $s_j$ replaced by $\max_{j \in G} s_j$. Let $\lambda_j \asymp \sqrt{\frac{\log p}{n}}$ uniformly in $j$ in (10). Assume that $\frac{(\log p)^7}{n} \leq C_1 n^{-c_1}$ for some constants $c_1, C_1 > 0$, $\max_{j=1,...,p} \frac{s_j^2(\log p)^3}{n} = o(1)$, $\frac{s_0^2(\log p)^3}{n} = o(1)$, and there exists a sequence of positive numbers $\alpha_n \to \infty$ such that $\alpha_n (\log p)^2 \max_{j=1,...,p} \lambda_j \sqrt{s_j} = o(1)$. Then under the null $H_{0,G}$, for any $G \subseteq \{1, 2, ..., p\}$, we have

$$\sup_{\alpha \in (0,1)} |\mathbb{P}(T_G > c_G(\alpha)) - \alpha| = o(1). \tag{14}$$

With Corollary 1, the power analysis of $T_G$ then follows from Theorem 2.4 in Zhang and Cheng (2017). The above testing procedure can be easily adapted for constructing simultaneous confidence intervals and support recovery, as we will see in Sections 5.2 and 5.3.

### 3.2 Theoretical implication of Theorem 1

The technique where we replace $X_{ij}$s by the estimated partial residuals $\hat{X}_{ij} = X_{ij} - \hat{\mathbb{E}}(X_{ij}|Z_i)$ as in (5)-(6) is called "partialling out". Note that this technique involves $p$ nonparametric regressions

where $p \geq n$. Moreover, the estimation error from each nonparametric regression accumulates in the approximate inverse $\hat{\Theta}$ of $n^{-1}\sum_{i=1}^{n} \hat{X}_i^T \hat{X}_i$. Consequently, we first discuss what makes the "partialling out" strategy work in the statistical inference of $\beta_0$ despite that $p$ is high dimensional.

Recall from our previous discussion that $\hat{b}_j - \beta_{0j}$ and $\tilde{b}_j - \beta_{0j}$ can be decomposed into a leading term $\frac{1}{n}\Theta_j \tilde{X}^T \varepsilon$ and several remainder terms as shown in (21)-(22). The rates of convergence for the remainder terms that are related to the nonparametric projection step depend on $\max_{j,j'} \left| \frac{1}{n}\sum_{i=1}^{n} \tilde{X}_{ij} \left[ \hat{f}_{j'}(Z_i) - f_{j'}(Z_i) \right] \right|$ with $\hat{f}_j$ defined in (7). In particular, we show that

$$\max_{j,j'} \left| \frac{1}{n}\sum_{i=1}^{n} \tilde{X}_{ij} \left[ \hat{f}_{j'}(Z_i) - f_{j'}(Z_i) \right] \right| \leq c t_n^2 \qquad \text{for any } t_n \geq r_n, \tag{15}$$

with probability at least $1 - \exp\left(-c' n t_n^2 + c'' \log p\right)$, for some constants $c, c', c'' > 0$. For many popular function classes, the *critical radius* $r_n$ defined earlier gives the optimal scaling for bounds on $\left\| \hat{f}_{j'} - f_{j'} \right\|_n$. In particular, for (7), one can show that

$$\max_{j'} \left\| \hat{f}_{j'} - f_{j'} \right\|_n \leq c'' t_n \qquad \text{for any } t_n \geq r_n, \tag{16}$$

with probability at least $1 - c'_0 \exp\left(-c'_1 n t_n^2 + c'_2 \log p\right)$.

Note that the orthogonality condition $\mathbb{E}\left(\tilde{X}_{ij}|Z_i\right) = 0$ (for all $j$) introduced by our partialling out strategy "reduces" the effects of the estimation errors from $\hat{f}_j$: The statistical error contributed by the projection step is $r_n^2$ instead of the optimal rate $r_n$ that one would expect from the nonparametric regression. Given this observation, for some function $h(s_j, s_0)$ of $s_j$ and $s_0$ only (where the exact form of $h$ is detailed in Assumption 5), as long as

$$\sqrt{n} r_n^2 h(s_j, s_0) = o(1),$$

the remainder terms related to (7) in the asymptotic expansions of $\sqrt{n}\left(\hat{b}_j - \beta_{0j}\right)$ and $\sqrt{n}\left(\tilde{b}_j - \beta_{0j}\right)$ are dominated by the leading term $\frac{1}{\sqrt{n}}\Theta_j \tilde{X}^T \varepsilon$, which has an asymptotic normal distribution. Note that the above finding also holds true for the surrogate $\hat{Y}_i = Y_i - \hat{\mathbb{E}}(Y_i|Z_i)$ (which is used in (5)) and the surrogate $\hat{g}(Z_i)$ (which is used in (6)).

### 3.3 Practical implication of Theorem 1

We illustrate the theoretical insight above with four specific examples in terms of $\dim(Z_i)$ and the function class $\mathcal{F}$ that $f_j$s belong to. The initial estimators $\hat{\beta}$ based on Zhu (2017) work for all four examples while $\hat{\beta}$ based on Müller and van de Geer (2015) or Yu, et al. (2017) only works for the first example. With regard to the specific forms of $\hat{\mathbb{E}}(X_{ij}|Z_i)$s, several modern techniques for the projection step are considered. The rates achieved by these practical procedures are then compared with the theoretical results in Section 3.2. To facilitate the presentation, our following discussions only concern $\mathbb{E}(X_{ij}|Z_i)$s and $\hat{\mathbb{E}}(X_{ij}|Z_i)$s; $\mathbb{E}(Y_i|Z_i)$ and $\hat{\mathbb{E}}(Y_i|Z_i)$ can be argued in the same fashion.

**Example 1**: $Z_i \in \mathbb{R}$ and $\mathcal{F} \in \mathcal{S}^m$ (the $m$th order Sobolev ball). Estimating $\mathbb{E}(X_{ij}|Z_i)$s via (7) or (8) can be reduced to the smoothing spline procedure, which achieves the sharp rate, $n^{-\frac{2m}{2m+1}}$, on $r_n^2$. In this case, we require $\sqrt{n} n^{-\frac{2m}{2m+1}} h(s_j, s_0) = o(1)$.

**Example 2**: $Z_i \in \mathbb{R}$ and $\mathcal{F}$ is the class of linear combinations of bounded basis functions $\psi_l(\cdot)$s such that for $f \in \mathcal{F}$, $f(Z_i) = \sum_{l=1}^{d_1} \theta_l \psi_l(Z_i)$ and $\theta := (\theta_l)_{l=1}^{d_1}$ belongs to the $l_0-$"ball" of "radius" $k$. Suppose $d_1 \geq n$ and $d_1 \geq 4k$. Then the standard Lasso procedure would yield upper bounds with scaling $\frac{k \log d_1}{n}$ on the quantities in (15) using the fact that

$$\left| \frac{1}{n} \sum_{i=1}^{n} \tilde{X}_{ij} \left[ \hat{f}_{j'}(Z_i) - f_{j'}(Z_i) \right] \right| \leq \max_{l=1,...,d_1} \left| \frac{1}{n} \sum_{i=1}^{n} \tilde{X}_{ij} \psi_l(Z_i) \right| \left\| \hat{\theta} - \theta \right\|_1$$

$$= O_p\left(\sqrt{\frac{\log d_1}{n}}\right) O_p\left(k\sqrt{\frac{\log d_1}{n}}\right)$$

$$= O_p\left(\frac{k \log d_1}{n}\right)$$

where $\hat{\theta}$ is the Lasso estimate. The scaling $\frac{k \log d_1}{n}$ almost achieves the sharp rate, $\frac{k \log \frac{d_1}{k}}{n}$, on $r_n^2$. In this case, we require $\frac{k \log d_1}{\sqrt{n}} h(s_j, s_0) = o(1)$.

If we use the recently proposed Slope (Su and Candès, 2016) instead of the standard Lasso, then the scaling $\frac{k \log \frac{d_1}{k}}{n}$ can be attained. In this case, we require $\frac{k \log \frac{d_1}{k}}{\sqrt{n}} h(s_j, s_0) = o(1)$.

**Example 3**: $Z_i \in \mathbb{R}$ and $\mathcal{F} \in \mathcal{RS}^m$ (the rearranged $m$th order Sobolev ball). Belloni, et al. (2014) show that, when used for estimating functions in an "enlarged" Sobolev space, i.e., $\mathcal{RS}^m$, the (square-root) Lasso achieves near oracle rates uniformly over the space. Hence, the square-root Lasso would be our estimation method in this example. The rearranged $m$th order Sobolev ball is defined as the class of functions which take the form $f(Z_i) = \sum_{l=1}^{\infty} \theta_l \psi_l(Z_i)$ such that $\sum_{l=1}^{\infty} |\theta_l| < \infty$ and

$$\theta := (\theta_1, \theta_2, ...) \in \Theta_R(m, d_1, L) = \left\{ \tilde{\theta} \in \ell^2(\mathbb{N}) : \begin{array}{c} \exists \text{ permutation } \Upsilon : \{1,...,d_1\} \to \{1,...,d_1\} \\ \sum_{l=1}^{d_1} l^{2m} \tilde{\theta}_{\Upsilon(l)}^2 + \sum_{l=d_1+1}^{\infty} l^{2m} \tilde{\theta}_l^2 \leq L \end{array} \right\}.$$

Note that $\mathcal{S}^m \subseteq \mathcal{RS}^m$. For technical purposes, define an oracle model $f^0 = \sum_{l=1}^{d_1} \theta_{0l} \psi_l(Z_i)$ in terms of a target parameter vector $\theta_0 := (\theta_{0l})_{l=1}^{d_1}$ that solves

$$\min_{\tilde{\theta} \in \mathbb{R}^{d_1}} \frac{1}{n} \sum_{i=1}^{n} \left[ f(Z_i) - \sum_{l=1}^{d_1} \tilde{\theta}_l \psi_l(Z_i) \right]^2 + \frac{\sigma^2 \left\| \tilde{\theta} \right\|_0}{n},$$

where $\left\| \tilde{\theta} \right\|_0$ denotes the number of non-zero components in $\tilde{\theta}$ and $\sigma^2 (\asymp 1)$ is the noise variance. The square-root Lasso $\hat{\theta}$ then provides an estimator for $\theta_0$ and the resulting nonparametric estimator for the true function $f(Z_i)$ is given by $\hat{f} = \sum_{l=1}^{d_1} \hat{\theta}_l \psi_l(Z_i)$; see Belloni, et al. (2014).

We employ the existing theoretical results to show

$$\left| \frac{1}{n} \sum_{i=1}^{n} \tilde{X}_{ij} \left[ \hat{f}_{j'}(Z_i) - f_{j'}(Z_i) \right] \right| \leq \left| \frac{1}{n} \sum_{i=1}^{n} \tilde{X}_{ij} \left[ \hat{f}_{j'}(Z_i) - f_{j'}^0(Z_i) \right] \right| + \left| \frac{1}{n} \sum_{i=1}^{n} \tilde{X}_{ij} \left[ f_{j'}^0(Z_i) - f_{j'}(Z_i) \right] \right|$$

$$\leq \max_{l=1,...,d_1} \left| \frac{1}{n} \sum_{i=1}^{n} \tilde{X}_{ij} \psi_l(Z_i) \right| \left\| \hat{\theta} - \theta_0 \right\|_1 + \left| \frac{1}{n} \sum_{i=1}^{n} \tilde{X}_{ij} \left[ f_{j'}^0(Z_i) - f_{j'}(Z_i) \right] \right|$$

$$= \underbrace{O_p\left(n^{\frac{1}{2m+1}} \frac{\log d_1}{n}\right)}_{T_1} + \underbrace{O_p\left(n^{-\frac{2m}{2m+1}}\right)}_{T_2} = O_p\left(n^{-\frac{2m}{2m+1}} \log d_1\right)$$

where $T_1$ follows from the bounds $\left\|\hat{\theta} - \theta_0\right\|_1 = O_p\left(n^{\frac{1}{2m+1}}\sqrt{\frac{\log d_1}{n}}\right)$ and $\max_l \left|\frac{1}{n}\sum_{i=1}^n \tilde{X}_{ij}\psi_l(Z_i)\right| = O_p\left(\sqrt{\frac{\log d_1}{n}}\right)$, $T_2$ follows from the bound $\left\|f_{j'} - f_{j'}^0\right\|_n = O_p\left(n^{-\frac{m}{2m+1}}\right)$ and (15). The bounds on $\left\|\hat{\theta} - \theta_0\right\|_1$ and $\left\|f_{j'} - f_{j'}^0\right\|_n$ are established in Belloni, et al. (2014). In this case, we require $\frac{1}{\sqrt{n}}\left(n^{-\frac{2m}{2m+1}}\log d_1\right)h(s_j, s_0) = o(1)$. Note that the scaling $n^{-\frac{2m}{2m+1}}\log d_1$ differs from the sharp rate, $n^{-\frac{2m}{2m+1}}$, on $r_n^2$ by a $\log d_1$ factor.

**Example 4**: $Z_i \in \mathbb{R}^d$ and $\mathcal{F}$ is the class of $|S| := k$ sparse additive nonparametric functions in the sense that any member $f$ in $\mathcal{F}$ has the following decomposition form: $f(Z_i) = \sum_{l=1}^d f_l(Z_{il}) = \sum_{l\in S} f_l(Z_{il})$; moreover, $f_l$ belongs to an RKHS of univariate functions. Suppose $d \geq n$ and $d \geq 4k$. In practice, we may apply the method in Raskutti, et al. (2012) to estimate $\mathbb{E}(X_{ij}|Z_i)$s. In particular, the estimators $\hat{f}_j$s defined in (7) can be written in the form

$$\hat{f}_j(z_1, ..., z_d) = \sum_{i=1}^n \sum_{l=1}^d \hat{\alpha}_{ilj}\mathbb{K}^l(z_l, z_{il})$$

where $\mathbb{K}^l$ denotes the kernel function for co-ordinate $l$ such that the collection of empirical kernel matrices $K^l \in \mathbb{R}^{n\times n}$ has entries $K^l_{ii'} = \mathbb{K}^l(z_{il}, z_{i'l})$; the optimal weights $\{\hat{\alpha}_{lj} \in \mathbb{R}^n, l = 1, ..., d\}$ are any solution to the following convex program proposed by Raskutti, et al. (2012):

$$\min_{\substack{\alpha_{lj} \in \mathbb{R}^n \\ \alpha_{lj}^T K^l \alpha_{lj} \leq 1}} \left\{\frac{1}{2n}\left\|w_j - \bar{w}_j - \sum_{l=1}^d K^l\alpha_{lj}\right\|_2^2 + \mu_{1j}\sum_{l=1}^d \sqrt{\frac{1}{n}\|K^l\alpha_{lj}\|_2^2} + \mu_{2j}\sum_{l=1}^d \sqrt{\alpha_{lj}^T K^l\alpha_{lj}}\right\}$$

where $\bar{w}_j = \frac{1}{n}\sum_{i=1}^n w_{ij}$ and $w_j = (w_{ij})_{i=1}^n$ (recalling from Section 2 that $w_{i0} = y_i$ and $w_{ij} = x_{ij}$ for each $j = 1, ..., p$).

If the underlying RKHS is $\mathcal{S}^m$, we would require $\sqrt{n}k\left(n^{-\frac{2m}{2m+1}} \vee \frac{\log d}{n}\right)h(s_j, s_0) = o(1)$ in this case. Note that the scaling $k\left(n^{-\frac{2m}{2m+1}} \vee \frac{\log d}{n}\right)$ almost achieves the sharp rate, $k\left(n^{-\frac{2m}{2m+1}} \vee \frac{\log \frac{d}{k}}{n}\right)$, on $r_n^2$.

### 3.4 The case of approximate sparsity

With additional efforts, the exact sparsity assumptions of $\beta_0$ and $\pi_j$ can be relaxed to accommodate for approximate sparsity, provided that the ordered coefficients decay sufficiently fast. To work with approximately sparse $\beta_0$ and $\pi_j$, we introduce two thresholded subsets:

$$S_{\underline{\tau}} := \{j \in \{1, 2, ..., p\} : |\beta_{0j}| > \underline{\tau}\}, \tag{17}$$

$$S_{\underline{\tau}_j} := \{l \in \{1, 2, ..., p\} \setminus j : |\pi_{jl}| > \underline{\tau}_j\}. \tag{18}$$

Let $|S_{\underline{\tau}}| := s_0$ and $|S_{\underline{\tau}_j}| := s_j$. Note that the newly defined $s_0$ and $s_j$ generalize the previous exact sparsity parameters; in Theorem 1, we simply take $\underline{\tau} = 0$ and $\underline{\tau}_j = 0$. Below we introduce Assumptions 4A and 5A. The roles these assumptions play in the case of approximately sparse $\beta_0$ and $\pi_j$ are similar to those Assumptions 4 and 5 play in the case of exact sparsity.

**Assumption 4A.** With $\underline{\tau} = \frac{c\sqrt{\frac{\log p}{n}}}{\Lambda_{\min}^2(\Sigma)}$ in (17) for some universal constant $c > 0$, (i) the initial estimator $\hat{\beta}$ satisfies that $\left\|\hat{\beta} - \beta_0\right\|_2 = O_p(D_2)$ and $\left\|\hat{\beta} - \beta_0\right\|_1 = O_p(D_1)$, where

$$D_1 = s_0\sqrt{\frac{\log p}{n}} + \left(s_0\sqrt{\frac{\log p}{n}}\left\|\beta_{0,S_{\underline{\tau}}^c}\right\|_1\right)^{\frac{1}{2}} + \left\|\beta_{0,S_{\underline{\tau}}^c}\right\|_1,$$

$$D_2 = \sqrt{\frac{s_0\log p}{n}} + \left(\sqrt{\frac{\log p}{n}}\left\|\beta_{0,S_{\underline{\tau}}^c}\right\|_1\right)^{\frac{1}{2}}.$$

(ii) $\left\|\mathbb{E}(X|Z)\left(\hat{\beta} - \beta_0\right)\right\|_n = O_p(D_2)$ and the estimator $\hat{g}$ satisfies that $\|\hat{g} - g\|_n^2 = O_p\left((\|\beta_0\|_1 r_n^2) \vee D_2^2\right)$.

**Remark.** Under mild conditions, the rates in Assumption 4A are satisfied by the initial estimators based on Zhu (2017). As Assumption 4(ii), Assumption 4A(ii) is only used in the analysis for the second debiased estimator $\tilde{b}$ in (6).

**Assumption 5A.** $(\|\pi_j\|_1 \vee 1) r_n^2 = O\left(\sqrt{\frac{\log p}{n}}\right)$ and

$$\sqrt{n}\max\left\{\|\pi_j\|_1 D_2^2,\ B_{1j}\sqrt{\frac{\log p}{n}},\ B_{1j}D_1,\ \|\pi_j\|_1(\|\beta_0\|_1 \vee 1)\left(r_n^2 \vee \frac{\log p}{n}\right)\right\} = o(1),$$

$$\sqrt{n}\max\left\{\left(\|\pi_j\|_1^2 \vee \|\pi_j\|_1\right)\frac{\log p}{n},\ \left(\|\pi_j\|_1^3 \vee \|\pi_j\|_1\right)\left(r_n^2 \vee \frac{\log p}{n}\right)\sqrt{\frac{\log p}{n}}\right\} = o(1),$$

$$\sqrt{n}\max\left\{\left(\|\pi_j\|_1^2 \vee \|\pi_j\|_1\right)\sqrt{\frac{\log p}{n}}D_1,\ \left(\|\pi_j\|_1^3 \vee \|\pi_j\|_1\right)\left(r_n^2 \vee \frac{\log p}{n}\right)D_1\right\} = o(1),$$

$$\max\left\{B_{1j},\ \left(\|\pi_j\|_1^2 \vee \|\pi_j\|_1\right)\sqrt{\frac{\log p}{n}},\ \left(\|\pi_j\|_1^3 \vee \|\pi_j\|_1\right)\left(r_n^2 \vee \frac{\log p}{n}\right)\right\} = o(1),$$

where

$$B_{1j} = s_j\sqrt{\frac{\log p}{n}} + \left(s_j\sqrt{\frac{\log p}{n}}\left\|\pi_{j,S_{\underline{\tau}_j}^c}\right\|_1\right)^{\frac{1}{2}} + \left\|\pi_{j,S_{\underline{\tau}_j}^c}\right\|_1,$$

$$B_{2j} = \sqrt{\frac{s_j \log p}{n}} + \left(\sqrt{\frac{\log p}{n}}\left\|\pi_{j,S_{\underline{\tau}_j}^c}\right\|_1\right)^{\frac{1}{2}}.$$

In comparison with the case of exact sparsity, approximate sparsity permits all the entries in $\beta_0$ (and $\pi_j$) to be non-zero at the expense of incurring an extra approximation error $\left(s_0\sqrt{\frac{\log p}{n}}\left\|\beta_{0,S_{\underline{\tau}}^c}\right\|_1\right)^{\frac{1}{2}} + \left\|\beta_{0,S_{\underline{\tau}}^c}\right\|_1$ (respectively, $\left(s_j\sqrt{\frac{\log p}{n}}\left\|\pi_{j,S_{\underline{\tau}_j}^c}\right\|_1\right)^{\frac{1}{2}} + \left\|\pi_{j,S_{\underline{\tau}_j}^c}\right\|_1$) in the upper bound for $\left\|\hat{\beta} - \beta_0\right\|_1$ (respectively, $\|\hat{\pi}_j - \pi_j\|_1$). Relative to Assumption 5, Assumption 5A imposes additional conditions on the

"small coefficients" via $\left\|\beta_{0,S_{\mathcal{T}}^c}\right\|_1$ and $\left\|\pi_{j,S_{\mathcal{T}_j}^c}\right\|_1$ so that the remainder terms in the asymptotic expansions of $\sqrt{n}\left(\hat{b}_j - \beta_{0j}\right)$ and $\sqrt{n}\left(\tilde{b}_j - \beta_{0j}\right)$ are dominated by the leading term $\frac{1}{\sqrt{n}}\Theta_j \tilde{X}^T \varepsilon$ (which gives Theorem 2 below). Note that for the special case where $\left\|\pi_{j,S_{\mathcal{T}_j}^c}\right\|_1 = O\left((s_j \vee 1)\sqrt{\frac{\log p}{n}}\right)$ and $\left\|\beta_{0,S_{\mathcal{T}}^c}\right\|_1 = O\left((s_0 \vee 1)\sqrt{\frac{\log p}{n}}\right)$, we have $D_1 = (s_0 \vee 1)\sqrt{\frac{\log p}{n}}$, $D_2 = \sqrt{\frac{(s_0 \vee 1)\log p}{n}}$, $B_{1j} = (s_j \vee 1)\sqrt{\frac{\log p}{n}}$, and $B_{2j} = \sqrt{\frac{(s_j \vee 1)\log p}{n}}$ in Assumptions 4A and 5A; these terms take the same forms as those used in Assumptions 4 and 5.

**Theorem 2.** *Under Assumptions 1-3, 4A, 5A and 6, if we choose $\lambda_j \asymp \sqrt{\frac{\log p}{n}}$ uniformly in $j$ in (10) and $\tau_j = \frac{c'\sqrt{\frac{\log p}{n}}}{\Lambda_{\min}^2(\Sigma)}$ in (18) for some universal constant $c' > 0$, then the claims in Theorem 1 still hold.*

## 4 Proofs

As a general rule for this appendix, all the $c \in (0, \infty)$ constants denote positive universal constants. The specific values of these constants may change from place to place. For notational simplicity, we assume the regime of interest is $p \geq (n \vee 2)$; the modification to allow $p < (n \vee 2)$ is trivial.

**Proof of Theorem 1 and Corollary 1**

Recall the two versions of the debiased estimators:

$$\hat{b}_j := \hat{\beta}_j + \frac{1}{n}\hat{\Theta}_j \hat{X}^T \left(\hat{Y} - \hat{X}\hat{\beta}\right), \tag{19}$$

$$\tilde{b}_j := \hat{\beta}_j + \frac{1}{n}\hat{\Theta}_j \hat{X}^T \left(Y - X\hat{\beta} - \hat{g}(Z)\right); \tag{20}$$

$\hat{g}(Z) := \{\hat{g}(Z_i)\}_{i=1}^n$, $\hat{Y} = Y - \hat{\mathbb{E}}(Y|Z) := \left\{Y_i - \hat{\mathbb{E}}(Y_i|Z_i)\right\}_{i=1}^n$ is an estimate for $\tilde{Y} = Y - \mathbb{E}(Y|Z) := \{Y_i - \mathbb{E}(Y_i|Z_i)\}_{i=1}^n$, and for $j = 1, ..., p$, $\hat{X}_j = X_j - \hat{\mathbb{E}}(X_j|Z) := \left\{X_{ij} - \hat{\mathbb{E}}(X_{ij}|Z_i)\right\}_{i=1}^n$ is an estimate for $\tilde{X}_j = X_j - \mathbb{E}(X_j|Z) = \{X_{ij} - \mathbb{E}(X_{ij}|Z_i)\}_{i=1}^n$. We write $\hat{Y} = \tilde{Y} + \hat{Y} - \tilde{Y} = \tilde{X}\beta_0 + \varepsilon + \hat{Y} - \tilde{Y}$ and $\hat{X} = \tilde{X} + \hat{X} - \tilde{X}$, which are used in the following derivations. We show in the following that $\hat{b}_j$ and $\tilde{b}_j$ have the same asymptotic distribution.

For (19), we obtain

$$\hat{b}_j - \beta_{0j}$$
$$= \hat{\beta}_j - \beta_{0j} + \frac{1}{n}\hat{\Theta}_j \hat{X}^T \left(\hat{Y} - \hat{X}\hat{\beta}\right)$$
$$= \frac{1}{n}\hat{\Theta}_j \tilde{X}^T \varepsilon - \frac{1}{n}\hat{\Theta}_j \left(\tilde{X} - \hat{X}\right)^T \varepsilon + \hat{\beta}_j - \beta_{0j} + \frac{1}{n}\hat{\Theta}_j \hat{X}^T \left[\tilde{X}\beta_0 - \left(\hat{X} - \tilde{X}\right)\beta_0 + \hat{Y} - \tilde{Y} - \hat{X}\hat{\beta}\right]$$
$$= \frac{1}{n}\Theta_j \tilde{X}^T \varepsilon + \underbrace{\frac{1}{n}\left(\hat{\Theta}_j - \Theta_j\right)\tilde{X}^T \varepsilon}_{E_0} - \underbrace{\frac{1}{n}\hat{\Theta}_j \left(\tilde{X} - \hat{X}\right)^T \varepsilon}_{E_1} + \underbrace{\left(e_j - \frac{1}{n}\hat{\Theta}_j \hat{X}^T \hat{X}\right)\left(\hat{\beta} - \beta_0\right)}_{E_2}$$
$$\underbrace{-\frac{1}{n}\hat{\Theta}_j \hat{X}^T \left(\hat{X} - \tilde{X}\right)\beta_0}_{E_3} + \underbrace{\frac{1}{n}\hat{\Theta}_j \hat{X}^T \left(\hat{Y} - \tilde{Y}\right)}_{E_4}. \tag{21}$$

For (20), we obtain

$$
\begin{aligned}
\tilde{b}_j - \beta_{0j} &= \hat{\beta}_j - \beta_{0j} + \frac{1}{n}\hat{\Theta}_j \hat{X}^T \left(Y - X\hat{\beta} - \hat{g}(Z)\right) \\
&= \frac{1}{n}\hat{\Theta}_j \tilde{X}^T \varepsilon - \frac{1}{n}\hat{\Theta}_j \left(\tilde{X} - \hat{X}\right)^T \varepsilon + \hat{\beta}_j - \beta_{0j} \\
&\quad + \frac{1}{n}\hat{\Theta}_j \hat{X}^T \left[\hat{X}\beta_0 - \left(\hat{X} - \tilde{X}\right)\beta_0 + (X - \tilde{X})\beta_0 - \hat{X}\hat{\beta} + \left(\hat{X} - \tilde{X}\right)\hat{\beta} - (X - \tilde{X})\hat{\beta} + g_0 - \hat{g}\right] \\
&= \frac{1}{n}\Theta_j \tilde{X}^T \varepsilon + \underbrace{\frac{1}{n}\left(\hat{\Theta}_j - \Theta_j\right)\tilde{X}^T \varepsilon}_{E_0} - \underbrace{\frac{1}{n}\hat{\Theta}_j \left(\tilde{X} - \hat{X}\right)^T \varepsilon}_{E_1} \\
&\quad + \underbrace{\left(e_j - \frac{1}{n}\hat{\Theta}_j \hat{X}^T \hat{X}\right)\left(\hat{\beta} - \beta_0\right)}_{E_2} + \underbrace{\frac{1}{n}\hat{\Theta}_j \hat{X}^T \left(\hat{X} - \tilde{X}\right)\left(\hat{\beta} - \beta_0\right)}_{E_3'} \\
&\quad - \underbrace{\frac{1}{n}\hat{\Theta}_j \tilde{X}^T \mathbb{E}(X|Z)\left(\hat{\beta} - \beta_0\right)}_{E_4'} - \underbrace{\frac{1}{n}\hat{\Theta}_j \left(\hat{X} - \tilde{X}\right)^T \mathbb{E}(X|Z)\left(\hat{\beta} - \beta_0\right)}_{E_5'} - \underbrace{\frac{1}{n}\hat{\Theta}_j \hat{X}^T (\hat{g} - g_0)}_{E_6'}. \quad (22)
\end{aligned}
$$

By elementary inequalities, we have

$$
\begin{aligned}
E_0 &\leq \left\|\hat{\Theta}_j - \Theta_j\right\|_1 \left\|\frac{1}{n}\tilde{X}^T \varepsilon\right\|_\infty, \\
E_1 &\leq \left\|\hat{\Theta}_j\right\|_1 \left\|\frac{1}{n}\left(\tilde{X} - \hat{X}\right)^T \varepsilon\right\|_\infty, \\
E_2 &\leq \left\|e_j - \frac{1}{n}\hat{\Theta}_j \hat{X}^T \hat{X}\right\|_\infty \left\|\hat{\beta} - \beta_0\right\|_1, \\
E_3 &\leq \left\|\hat{\Theta}_j\right\|_1 \left\|\frac{1}{n}\hat{X}^T \left(\hat{X} - \tilde{X}\right)\right\|_\infty \|\beta_0\|_1, \\
E_4 &\leq \left\|\hat{\Theta}_j\right\|_1 \left\|\frac{1}{n}\hat{X}^T \left(\hat{Y} - \tilde{Y}\right)\right\|_\infty,
\end{aligned}
$$

and

$$
\begin{aligned}
E_3' &\leq \left\|\hat{\Theta}_j\right\|_1 \left\|\frac{1}{n}\hat{X}^T \left(\hat{X} - \tilde{X}\right)\right\|_\infty \left\|\hat{\beta} - \beta_0\right\|_1, \\
E_4' &\leq \left\|\hat{\Theta}_j\right\|_1 \left\|\frac{1}{n}\tilde{X}^T \mathbb{E}(X|Z)\right\|_\infty \left\|\hat{\beta} - \beta_0\right\|_1, \\
E_5' &\leq \left\|\hat{\Theta}_j\right\|_1 \max_{j=1,\ldots,p} \left\|\hat{X}_j - \tilde{X}_j\right\|_n \left\|\mathbb{E}(X|Z)\left(\hat{\beta} - \beta_0\right)\right\|_n, \\
E_6' &\leq \left\|\hat{\Theta}_j\right\|_1 \left\|\frac{1}{n}\hat{X}^T (\hat{g} - g_0)\right\|_\infty.
\end{aligned}
$$

We bound $\left\|e_j - \frac{1}{n}\hat{\Theta}_j \hat{X}^T \hat{X}\right\|_\infty$ with (41) and $\left\|\hat{\Theta}_j - \Theta_j\right\|_1$ with (38), which also implies that

$$
\left\|\hat{\Theta}_j\right\|_1 = O_p\left(\max_j s_j\right) + O_p\left(\max_j s_j^2 \sqrt{\frac{\log p}{n}} + \max_j s_j^3 \left(r_{nj}^2 \vee \frac{\log p}{n}\right)\right). \quad (23)
$$

By Assumption 4, we have

$$\left\|\hat{\beta} - \beta_0\right\|_1 = O_p\left(s_0\sqrt{\frac{\log p}{n}}\right),$$

$$\left\|\mathbb{E}(X|Z)\left(\hat{\beta} - \beta_0\right)\right\|_n = O_p\left(\sqrt{\frac{s_0 \log p}{n}}\right). \tag{24}$$

By Assumption 1, standard tail bounds for sub-Exponential variables [e.g., Vershynin (2012)] yield

$$\left\|\frac{1}{n}\tilde{X}^T \varepsilon\right\|_\infty = O_p\left(\sqrt{\frac{\log p}{n}}\right),$$

$$\left\|\frac{1}{n}\tilde{X}^T \mathbb{E}(X|Z)\right\|_\infty = O_p\left(\sqrt{\frac{\log p}{n}}\right), \tag{25}$$

where we have used the fact that $\mathbb{E}(X_{ij}|Z_i) = X_{ij} - \tilde{X}_{ij}$ is sub-Gaussian [implied by Assumption 1 and that sub-Gaussianity is preserved by linear operations]. Note that (24)-(25) only matter to the second debiased estimator $\tilde{b}_j$ in (20).

Note that we have

$$\left\|\hat{X}_j - \tilde{X}_j\right\|_n = \left\|\mathbb{E}(X_j|Z) - \hat{\mathbb{E}}(X_j|Z)\right\|_n,$$

$$\left|\frac{1}{n}\left(\tilde{X}_j - \hat{X}_j\right)^T \varepsilon\right| = \left|\frac{1}{n}\left[\hat{\mathbb{E}}(X_j|Z) - \mathbb{E}(X_j|Z)\right]^T \varepsilon\right|,$$

$$\left|\frac{1}{n}\hat{X}_j^T \left(\hat{X}_{j'} - \tilde{X}_{j'}\right)\right| \leq \left\|\mathbb{E}(X_j|Z) - \hat{\mathbb{E}}(X_j|Z)\right\|_n \left\|\mathbb{E}(X_{j'}|Z) - \hat{\mathbb{E}}(X_{j'}|Z)\right\|_n$$
$$+ \left|\frac{1}{n}\tilde{X}_j^T \left[\hat{\mathbb{E}}(X_{j'}|Z) - \mathbb{E}(X_{j'}|Z)\right]\right|,$$

$$\left|\frac{1}{n}\hat{X}_j^T \left(\hat{Y} - \tilde{Y}\right)\right| \leq \left\|\mathbb{E}(X_j|Z) - \hat{\mathbb{E}}(X_j|Z)\right\|_n \left\|\mathbb{E}(Y|Z) - \hat{\mathbb{E}}(Y|Z)\right\|_n + \left|\frac{1}{n}\tilde{X}_j^T \left[\hat{\mathbb{E}}(Y|Z) - \mathbb{E}(Y|Z)\right]\right|,$$

$$\left|\frac{1}{n}\hat{X}_j^T \left(\hat{g} - g_0\right)\right| \leq \left\|\mathbb{E}(X_j|Z) - \hat{\mathbb{E}}(X_j|Z)\right\|_n \|\hat{g} - g_0\|_n + \left|\frac{1}{n}\tilde{X}_j^T \left(\hat{g} - g_0\right)\right|.$$

We write $\hat{f}_j := \hat{\mathbb{E}}(X_j|Z)$ and $f_j := \mathbb{E}(X_j|Z)$ for $j = 1, ..., p$, as well as $\hat{f}_0 := \hat{\mathbb{E}}(Y|Z)$ and $f_0 := \mathbb{E}(Y|Z)$. Under Assumptions 1 and 3, by standard argument for nonparametric least squares estimators, for any $t_n \geq r_n$,

$$\left\|\hat{f}_j - f_j\right\|_n \leq c_1 t_n \tag{26}$$

with probability at least $1 - c\exp\left(-c' n t_n^2\right)$. Moreover, under Assumption 4,

$$\|\hat{g} - g_0\|_n = O_p\left(\sqrt{s_0} r_n \vee \sqrt{\frac{s_0 \log p}{n}}\right). \tag{27}$$

For fixed elements $\tilde{f}_j \in \mathcal{F}$ and $\tilde{g} \in \mathcal{F}$, respectively, we can view $\tilde{f}_j - f_j$ and $\tilde{g} - g_0$ as functions of $Z$ only. Additionally, note that $\mathbb{E}(\varepsilon_i|Z_i) = 0$ and $\mathbb{E}\left(\tilde{Y}_i|Z_i\right) = \mathbb{E}\left(\tilde{X}_{ij}|Z_i\right) = 0$ (by construction of $\tilde{Y}_i$ and $\tilde{X}_{ij}$, $i = 1, ..., n$, $j = 1, ..., p$). The remaining argument uses results from empirical processes theory and local function complexity. In particular, under the independent sampling assumption, Lemma A5 with (26)-(27) implies that

$$\max_{j=1,\ldots,p} \left| \frac{1}{n} \left[ \mathbb{E}(X_j|Z) - \hat{\mathbb{E}}(X_j|Z) \right]^T \varepsilon \right| = O_p(t_n^2),$$

$$\max_{j,j' \in \{1,\ldots,p\}} \left| \frac{1}{n} \tilde{X}_j^T \left[ \hat{\mathbb{E}}(X_{j'}|Z) - \mathbb{E}(X_{j'}|Z) \right] \right| = O_p(t_n^2),$$

$$\max_{j=1,\ldots,p} \left| \frac{1}{n} \tilde{X}_j^T \left[ \hat{\mathbb{E}}(Y|Z) - \mathbb{E}(Y|Z) \right] \right| = O_p(t_n^2),$$

$$\max_{j=1,\ldots,p} \left| \frac{1}{n} \tilde{X}_j^T (\hat{g} - g_0) \right| = O_p(s_0 t_n^2) + O_p\left(\frac{s_0 \log p}{n}\right),$$

provided that $t_n \geq r_n$ and $nt_n^2 \gtrsim \log p$. In our analysis, it the suffices to choose $t_n = \left(r_n \vee \sqrt{\frac{\log p}{n}}\right)$. Consequently,

$$\left\| \frac{1}{n} \left( \tilde{X} - \hat{X} \right)^T \varepsilon \right\|_\infty = O_p\left(r_n^2 \vee \frac{\log p}{n}\right),$$

$$\left\| \frac{1}{n} \hat{X}^T \left( \hat{X} - \tilde{X} \right) \right\|_\infty = O_p\left(r_n^2 \vee \frac{\log p}{n}\right), \tag{28}$$

$$\left\| \frac{1}{n} \hat{X}^T \left( \hat{Y} - \tilde{Y} \right) \right\|_\infty = O_p\left(r_n^2 \vee \frac{\log p}{n}\right), \tag{29}$$

$$\left\| \frac{1}{n} \hat{X}^T (\hat{g} - g_0) \right\|_\infty = O_p\left(s_0 r_n^2 \vee \frac{s_0 \log p}{n}\right). \tag{30}$$

Note that (29) only matters to the first debiased estimator $\hat{b}_j$ in (19) while (30) only matters to the second debiased estimator $\tilde{b}_j$ in (20).

Putting all the pieces together, under Assumption 5, we apply the CLT and obtain $\frac{\sqrt{n}(\hat{b}_j - \beta_{0j})}{\sigma_j} \xrightarrow{D} \mathcal{N}(0,1)$ and $\frac{\sqrt{n}(\tilde{b}_j - \beta_{0j})}{\sigma_j} \xrightarrow{D} \mathcal{N}(0,1)$, where $\sigma_j^2 = \Theta_j \mathbb{E}\left(\tilde{X}_i^T \tilde{X}_i\right) \Theta_j^T$. Now it remains to show that

$$\left| \hat{\Theta}_j \frac{\hat{X}^T \hat{X}}{n} \hat{\Theta}_j^T - \Theta_j \mathbb{E}\left(\tilde{X}_i^T \tilde{X}_i\right) \Theta_j^T \right|$$

$$\leq \left| \hat{\Theta}_j \left[ \frac{\hat{X}^T \hat{X}}{n} - \mathbb{E}\left(\tilde{X}_i^T \tilde{X}_i\right) \right] \hat{\Theta}_j^T \right| + \left| \hat{\Theta}_j \mathbb{E}\left(\tilde{X}_i^T \tilde{X}_i\right) \hat{\Theta}_j^T - \Theta_j \mathbb{E}\left(\tilde{X}_i^T \tilde{X}_i\right) \Theta_j^T \right|$$

$$\leq \left\| \frac{\hat{X}^T \hat{X}}{n} - \mathbb{E}\left(\tilde{X}_i^T \tilde{X}_i\right) \right\|_\infty \left\| \hat{\Theta} \right\|_1^2 + \left| \hat{\Theta}_j \mathbb{E}\left(\tilde{X}_i^T \tilde{X}_i\right) \hat{\Theta}_j^T - \Theta_j \mathbb{E}\left(\tilde{X}_i^T \tilde{X}_i\right) \Theta_j^T \right|$$

$$= O_p\left( s_j^2 \left( r_n^2 \vee \frac{\log p}{n} \right) + s_j^2 \sqrt{\frac{\log p}{n}} + s_j^3 \left( r_{nj}^2 \vee \frac{\log p}{n} \right) \right) \tag{31}$$

$$+ O_p\left( s_j^3 \frac{\log p}{n} + s_j^5 \left( r_n^2 \vee \frac{\log p}{n} \right)^2 + (s_j \vee 1) \sqrt{\frac{\log p}{n}} + \left( s_j^2 \vee 1 \right) \left( r_n^2 \vee \frac{\log p}{n} \right) \right) \tag{32}$$

where (31) follows from (23) and (47)-(48); (32) follows from (40). Thus we have shown Theorem 1.

To show Corollary 1, we adopt the following result (Lemma A1) from Lemma 1.1 of Zhang and Cheng (2016). Combining Lemma A1 with the facts established in Theorem 1, the claim in Corollary 1 follows.

**Lemma A1**. Let $\zeta_j = (\zeta_{1j}, ..., \zeta_{nj})$, $j = 1, ..., p$, be centered i.i.d. sub-Gaussian random variables with variance $\Theta_j \mathbb{E}\left(\tilde{X}_i^T \tilde{X}_i\right) \Theta_j^T$. Assume that $\frac{(\log p)^7}{n} \leq C_1 n^{-c_1}$ for some constants $c_1, C_1 > 0$; $\max_{j=1,...,p} \frac{s_j^2 (\log p)^3}{n} = o(1)$; $\frac{s_0^2 (\log p)^3}{n} = o(1)$; and there exists a sequence of positive numbers $\alpha_n \to \infty$ such that $\alpha_n (\log p)^2 \max_{j=1,...,p} \lambda_j \sqrt{s_j} = o(1)$. Then, for any $G \subseteq \{1, ..., p\}$,

$$\sup_{t \in \mathbb{R}} \left| \mathbb{P}\left(\frac{\max_{j \in G} \zeta_{ij}}{\sqrt{n}} \leq t\right) - \mathbb{P}\left(\frac{\max_{j \in G} \epsilon_{ij}}{\sqrt{n}} \leq t\right) \right| \lesssim n^{-c'}, \quad c' > 0,$$

where $\epsilon_j = (\epsilon_{j1}, ..., \epsilon_{jn})$ are i.i.d. normal random variables with mean 0 and variance $\Theta_j \mathbb{E}\left(\tilde{X}_i^T \tilde{X}_i\right) \Theta_j^T$.

**Lemma A2**. Suppose Assumptions 1, 2(i) regarding $\Lambda^2_{\min}(\Sigma)$, 3, and 6 hold. If $(\|\pi_j\|_1 \vee 1) r_n^2 = O\left(\sqrt{\frac{\log p}{n}}\right)$, $\lambda_j \gtrsim \sqrt{\frac{\log p}{n}}$, and

$$\max_j s_j \left[r_n^2 \vee \frac{\log p}{n}\right] \leq c \Lambda^2_{\min}(\Sigma) \tag{33}$$

for some sufficiently small constant $c > 0$, then,

$$\max_j \|\hat{\pi}_j - \pi_j\|_2 = O_p\left(\sqrt{s_j} \lambda_j\right), \tag{34}$$

$$\max_j \|\hat{\pi}_j - \pi_j\|_1 = O_p\left(s_j \lambda_j\right). \tag{35}$$

**Proof**. First, write $\eta_j := \tilde{X}_j - \tilde{X}_{-j} \pi_j$ and

$$\begin{aligned}
\tilde{X}_j &= \hat{X}_j + \hat{\mathbb{E}}(X_j | Z) - \mathbb{E}(X_j | Z) \\
&= \tilde{X}_{-j} \pi_j + \eta_j = \left[\hat{X}_{-j} + \hat{\mathbb{E}}(X_{-j} | Z) - \mathbb{E}(X_{-j} | Z)\right] \pi_j + \eta_j,
\end{aligned}$$

thus we have

$$\begin{aligned}
\hat{X}_j &= \hat{X}_{-j} \pi_j + \left[\hat{\mathbb{E}}(X_{-j} | Z) - \mathbb{E}(X_{-j} | Z)\right] \pi_j - \left[\hat{\mathbb{E}}(X_j | Z) - \mathbb{E}(X_j | Z)\right] + \eta_j \\
&= \hat{X}_{-j} \pi_j + u_j.
\end{aligned}$$

where

$$u_j := \left[\hat{\mathbb{E}}(X_{-j} | Z) - \mathbb{E}(X_{-j} | Z)\right] \pi_j - \left[\hat{\mathbb{E}}(X_j | Z) - \mathbb{E}(X_j | Z)\right] + \eta_j. \tag{36}$$

By standard argument for the Lasso, applying Lemma A6 [which shows the $\frac{1}{n} \hat{X}_{-j}^T \hat{X}_{-j}$ satisfies a lower restricted eigenvalue (LRE) condition with probability at least $1 - o(1)$] and Lemma A7 along with Assumption 6 [which implies that $\max_j \left\|\frac{1}{n} \hat{X}_{-j}^T u_j\right\|_\infty = O_p\left(\sqrt{\frac{\log p}{n}}\right)$] yields (34) and (35).

**Lemma A3**. Suppose Assumptions in Lemma A2 hold. Let $\lambda_j \asymp \sqrt{\frac{\log p}{n}}$ uniformly in $j$. Then for every $j = 1, ..., p$, we have

$$\left|\hat{\tau}_j^2 - \tau_j^2\right| = O_p\left(\max_j (s_j \vee 1) \sqrt{\frac{\log p}{n}} + \max_j (s_j^2 \vee 1) \left(r_n^2 \vee \frac{\log p}{n}\right)\right), \tag{37}$$

where $\tau_j^2 := \mathbb{E}\left[\left(\tilde{X}_{ij} - \tilde{X}_{i,-j}\pi_j\right)^2\right]$; moreover,

$$\left\|\hat{\Theta}_j - \Theta_j\right\|_1 = O_p\left(\max_j s_j^2 \sqrt{\frac{\log p}{n}} + \max_j s_j^3 \left(r_n^2 \vee \frac{\log p}{n}\right)\right), \tag{38}$$

$$\left\|\hat{\Theta}_j - \Theta_j\right\|_2 = O_p\left(\max_j s_j^{\frac{3}{2}} \sqrt{\frac{\log p}{n}} + \max_j s_j^{\frac{5}{2}} \left(r_n^2 \vee \frac{\log p}{n}\right)\right), \tag{39}$$

$$\left|\hat{\Theta}_j \mathbb{E}\left(\tilde{X}_i^T \tilde{X}_i\right) \hat{\Theta}_j^T - \Theta_{j,j}\right| = O_p\left(\max_j s_j^3 \frac{\log p}{n} + \max_j s_j^5 \left(r_n^2 \vee \frac{\log p}{n}\right)^2\right)$$

$$+ O_p\left(\max_j (s_j \vee 1) \sqrt{\frac{\log p}{n}} + \max_j (s_j^2 \vee 1)\left(r_n^2 \vee \frac{\log p}{n}\right)\right), \tag{40}$$

$$\left\|\hat{\Theta}_j \frac{\hat{X}^T \hat{X}}{n} - e_j\right\|_\infty = O_p\left(\max_j s_j^2 \sqrt{\frac{\log p}{n}} + \max_j s_j^3 \left(r_n^2 \vee \frac{\log p}{n}\right)\right). \tag{41}$$

**Proof.** Note that we have $\hat{\tau}_j^2 := \frac{1}{n}\left\|\hat{X}_j - \hat{X}_{-j}\hat{\pi}_j\right\|_2^2 + \lambda_j \left\|\hat{\pi}_j\right\|_1$ and

$$\left|\frac{1}{n}\sum_{i=1}^n \left(\hat{X}_{ij} - \hat{X}_{i,-j}\hat{\pi}_j\right)^2 - \tau_j^2\right|$$

$$\leq \frac{1}{n}\sum_{i=1}^n \left[\hat{X}_{i,-j}(\pi_j - \hat{\pi}_j)\right]^2 + \left|\frac{2}{n}\sum_{i=1}^n \left[u_{ij}\hat{X}_{i,-j}(\pi_j - \hat{\pi}_j)\right]\right| + \left|\frac{1}{n}\sum_{i=1}^n u_{ij}^2 - \tau_j^2\right|$$

$$\leq \frac{2}{n}\sum_{i=1}^n \left[\tilde{X}_{i,-j}(\pi_j - \hat{\pi}_j)\right]^2 + \|\pi_j - \hat{\pi}_j\|_1^2 \left\|\frac{2}{n}\sum_{i=1}^n \left(\hat{X}_{i,-j} - \tilde{X}_{i,-j}\right)^T \left(\hat{X}_{i,-j} - \tilde{X}_{i,-j}\right)\right\|_\infty$$

$$+ \left\|\frac{2}{n}\hat{X}_{-j}^T u_j\right\|_\infty \|\pi_j - \hat{\pi}_j\|_1 + \left|\frac{1}{n}\sum_{i=1}^n \left(u_{ij}^2 - \eta_{ij}^2\right)\right| + \left|\frac{1}{n}\sum_{i=1}^n \left(\eta_{ij}^2 - \tau_j^2\right)\right|.$$

Under the condition $\Lambda_{\max}^2(\Sigma) = O(1)$ [Assumption 2], applying (34) in Lemma A2 and Lemma A4 yields

$$\frac{2}{n}\sum_{i=1}^n \left[\tilde{X}_{i,-j}(\pi_j - \hat{\pi}_j)\right]^2 = O_p\left(\frac{\max_j s_j \log p}{n}\right).$$

By choosing $t_n = \left(r_n \vee \sqrt{\frac{\log p}{n}}\right)$ as in the proof for Theorem 1, (26) and (35) imply that

$$\|\pi_j - \hat{\pi}_j\|_1^2 \left\|\frac{2}{n}\sum_{i=1}^n \left(\hat{X}_{i,-j} - \tilde{X}_{i,-j}\right)^T \left(\hat{X}_{i,-j} - \tilde{X}_{i,-j}\right)\right\|_\infty$$

$$= O_p\left(\max_j s_j^2 \frac{\log p}{n}\right) O_p\left(r_n^2 \vee \frac{\log p}{n}\right).$$

By (35) and Lemma A7 along with Assumption 6, we have

$$\left\|\frac{2}{n}\hat{X}_{-j}^T u_j\right\|_\infty \|\pi_j - \hat{\pi}_j\|_1 = O_p\left(\sqrt{\frac{\log p}{n}}\right) O_p\left(\max_j s_j \sqrt{\frac{\log p}{n}}\right).$$

For the term $\left|\frac{1}{n}\sum_{i=1}^{n}\left(u_{ij}^2 - \eta_{ij}^2\right)\right|$, it suffices to show

$$\left\|\left[\hat{\mathbb{E}}\left(X_{-j}|Z\right) - \mathbb{E}\left(X_{-j}|Z\right)\right]\pi_j\right\|_n^2 = O_p\left(\max_j s_j^2\left(r_n^2 \vee \frac{\log p}{n}\right)\right) \tag{42}$$

$$\left\|\hat{\mathbb{E}}\left(X_j|Z\right) - \mathbb{E}\left(X_j|Z\right)\right\|_n^2 = O_p\left(r_n^2 \vee \frac{\log p}{n}\right) \tag{43}$$

$$\frac{1}{n}\eta_j^T\left[\hat{\mathbb{E}}\left(X_{-j}|Z\right) - \mathbb{E}\left(X_{-j}|Z\right)\right]\pi_j = O_p\left(\max_j s_j\left(r_n^2 \vee \frac{\log p}{n}\right)\right) \tag{44}$$

$$\frac{1}{n}\eta_j^T\left[\hat{\mathbb{E}}\left(X_j|Z\right) - \mathbb{E}\left(X_j|Z\right)\right] = O_p\left(r_n^2 \vee \frac{\log p}{n}\right). \tag{45}$$

In the above, (42) and (43) follow from (26) (where we choose $t_n = \left(r_n \vee \sqrt{\frac{\log p}{n}}\right)$). In terms of (44) and (45), for fixed elements $\tilde{f}_{j'}, f_{j'} \in \mathcal{F}$, we can view $\tilde{f}_{j'} - f_{j'}$ as functions of $Z$ only. Meanwhile, $\eta_j$ is a function of $\tilde{X}$ only, so $\mathbb{E}\left(\tilde{X}_{ij}|Z_i\right) = 0$ ($i = 1, ..., n$, $j = 1, ..., p$) implies that

$$\mathbb{E}\left[\eta_{ij}\tilde{f}_{j'}(Z_i)|Z_i\right] = \tilde{f}_{j'}(Z_i)\mathbb{E}\left[\tilde{X}_{ij}|Z_i\right] - \tilde{f}_{j'}(Z_i)\mathbb{E}\left[\tilde{X}_{i,-j}|Z_i\right]\pi_j = 0. \tag{46}$$

Furthermore, $\mathbb{E}\left[\eta_{ij}\tilde{f}_{j'}(Z_i)\right] = 0$. As for (28), we apply Lemma A5 with (26) and obtain the (44) and (45).

In terms of $\left|\frac{1}{n}\sum_{i=1}^{n}\left(\eta_{ij}^2 - \tau_j^2\right)\right|$, we note that by Assumption 1 and the definition of $\eta_{ij}$, for $j = 1, ..., p$, an application of standard tail bounds for sub-Exponential variables yields

$$\left|\frac{1}{n}\sum_{i=1}^{n}\left(\eta_{ij}^2 - \tau_j^2\right)\right| = O_p\left(\sqrt{\frac{\log p}{n}}\right).$$

Moreover, by (35) and the choice $\lambda_j \asymp \sqrt{\frac{\log p}{n}}$, for $j = 1, ..., p$, we have

$$\lambda_j \|\hat{\pi}_j\|_1 = O_p\left(\max_j s_j \sqrt{\frac{\log p}{n}}\right) + O_p\left(\sqrt{\frac{\log p}{n}}\right) O\left(\max_j s_j \sqrt{\frac{\log p}{n}}\right)$$

Putting the pieces together, we have (37).

Next we show (38) and (39). Note that

$$\left\|\hat{\Theta}_j - \Theta_j\right\|_1 = \left\|\frac{\hat{M}_j}{\hat{\tau}_j^2} - \frac{M_j}{\tau_j^2}\right\|_1 \leq \frac{\|\hat{\pi}_j - \pi_j\|_1}{\hat{\tau}_j^2} + \|\pi_j\|_1\left(\frac{1}{\hat{\tau}_j^2} - \frac{1}{\tau_j^2}\right)$$

where the first term is $O_p\left(s_j\sqrt{\frac{\log p}{n}}\right)$ by (35) and the fact that $\hat{\tau}_j^2 = \tau_j^2 + o_p(1)$ while $\frac{1}{\tau_j^2} \lesssim 1$ [by Assumption 2]. For the second term, we have $\|\pi_j\|_1 = O(s_j)$ and

$$\frac{1}{\hat{\tau}_j^2} - \frac{1}{\tau_j^2} = O_p\left\{\left[\max_j(s_j \vee 1)\sqrt{\frac{\log p}{n}}\right] \vee \left[\max_j(s_j^2 \vee 1)\left(r_n^2 \vee \frac{\log p}{n}\right)\right]\right\}.$$

Therefore, we have (38). Similarly, we can obtain (39) by exploiting

$$\left\|\hat{\Theta}_j - \Theta_j\right\|_2 \leq \frac{\|\hat{\pi}_j - \pi_j\|_2}{\hat{\tau}_j^2} + \|\pi_j\|_2\left(\frac{1}{\hat{\tau}_j^2} - \frac{1}{\tau_j^2}\right).$$

We now show (41). Note that

$$\left\|\hat{\Theta}_j \frac{\hat{X}^T \hat{X}}{n} - e_j\right\|_\infty \leq \left\|\hat{\Theta}_j - \Theta_j\right\|_1 \left\|\frac{\hat{X}^T \hat{X}}{n}\right\|_\infty + \|\Theta_j\|_1 \left\|\frac{\hat{X}^T \hat{X}}{n} - \mathbb{E}\left(\tilde{X}_i^T \tilde{X}_i\right)\right\|_\infty.$$

By (28), we have

$$\left\|\frac{\hat{X}^T \hat{X}}{n} - \frac{\tilde{X}^T \tilde{X}}{n}\right\|_\infty \leq \left\|\frac{\tilde{X}^T (\hat{X} - \tilde{X})}{n}\right\|_\infty + \left\|\frac{(\hat{X} - \tilde{X})^T \tilde{X}}{n}\right\|_\infty + \left\|\frac{(\hat{X} - \tilde{X})^T (\hat{X} - \tilde{X})}{n}\right\|_\infty$$

$$= O_p\left(r_n^2 \vee \frac{\log p}{n}\right). \tag{47}$$

Moreover, by standard tail bounds for sub-Exponential variables, we have

$$\left\|\frac{\tilde{X}^T \tilde{X}}{n} - \mathbb{E}\left(\tilde{X}_i^T \tilde{X}_i\right)\right\|_\infty = O_p\left(\sqrt{\frac{\log p}{n}}\right). \tag{48}$$

Consequently, we obtain (41).

Finally we show (40). Using the facts that $\Theta_j \mathbb{E}\left(\tilde{X}_i^T \tilde{X}_i\right) = e_j^T$, $\Theta_j \mathbb{E}\left(\tilde{X}_i^T \tilde{X}_i\right) \Theta_j^T = \Theta_{j,j}$, $\hat{\Theta}_{j,j} = \frac{1}{\hat{\tau}_j^2}$, and $\Theta_{j,j} = \frac{1}{\tau_j^2}$, we have

$$\hat{\Theta}_j \mathbb{E}\left(\tilde{X}_i^T \tilde{X}_i\right) \hat{\Theta}_j^T - \Theta_{j,j}$$
$$= (\hat{\Theta}_j - \Theta_j) \mathbb{E}\left(\tilde{X}_i^T \tilde{X}_i\right) (\hat{\Theta}_j - \Theta_j)^T + 2\Theta_j \mathbb{E}\left(\tilde{X}_i^T \tilde{X}_i\right) (\hat{\Theta}_j - \Theta_j)^T + \Theta_j \mathbb{E}\left(\tilde{X}_i^T \tilde{X}_i\right) \Theta_j^T - \Theta_{j,j}$$
$$= (\hat{\Theta}_j - \Theta_j) \mathbb{E}\left(\tilde{X}_i^T \tilde{X}_i\right) (\hat{\Theta}_j - \Theta_j)^T + \frac{2}{\hat{\tau}_j^2} - \frac{2}{\tau_j^2}.$$

Note that

$$(\hat{\Theta}_j - \Theta_j) \mathbb{E}\left(\tilde{X}_i^T \tilde{X}_i\right) (\hat{\Theta}_j - \Theta_j)^T \leq \Lambda_{\max}^2 \left\|\hat{\Theta}_j - \Theta_j\right\|_2^2.$$

Applying (39) yields the claim. □

**Lemma A4.** Let $U \in \mathbb{R}^{n \times p_1}$ be a sub-Gaussian matrix with parameter $\sigma_U$ and each row be independently sampled. Then

$$\frac{\|U\Delta\|_2^2}{n} \leq \frac{3\bar{\alpha}}{2} \|\Delta\|_2^2 + \alpha' \frac{\log p_1}{n} \|\Delta\|_1^2, \quad \text{for all } \Delta \in \mathbb{R}^{p_1}$$

$$\frac{\|U\Delta\|_2^2}{n} \geq \frac{\underline{\alpha}}{2} \|\Delta\|_2^2 - \alpha' \frac{\log p_1}{n} \|\Delta\|_1^2, \quad \text{for all } \Delta \in \mathbb{R}^{p_1}$$

with probability at least $1 - c_1 \exp(-bn)$, where $\underline{\alpha}$, $\bar{\alpha}$, $\alpha'$, and $b$ are positive constants that only depend on $\sigma_U$ and the eigenvalues of $\Sigma_U = \mathbb{E}(U_i^T U_i)$ for $i = 1, ..., n$.

**Remark.** This lemma is Lemma 13 in Loh and Wainwright (2012).

**Lemma A5.** Suppose $\{v_i\}_{i=1}^n$ are i.i.d. sub-Gaussian variables with parameter at most $\sigma_v$ and $\mathbb{E}(v_i | Z_i) = 0$. Let $\bar{\mathcal{F}}$ be a star-shaped function class. Then, there are universal positive constants $(c, c_1, c_2)$ such that for any $t_n \geq r_n$ (where $r_n$ is the *critical radius*),

$$\sup_{f \in \bar{\mathcal{F}}, \|f\|_n \leq t_n} \left|\frac{1}{n} \sum_{i=1}^n v_i f(Z_i)\right| \leq c t_n^2 \tag{49}$$

with probability at least $1 - c_1 \exp\left(-c_2 n \frac{t_n^2}{\sigma_v^2}\right)$. If $\bar{\mathcal{F}}$ concerns a ball of a RKHS $\mathcal{H}$ equipped with a norm $\|\cdot\|_\mathcal{H}$, then for any $f \in \bar{\mathcal{F}}$ and $t_n \geq r_n$, we have

$$\left|\frac{1}{n}\sum_{i=1}^n v_i f(Z_i)\right| \leq c' t_n^2 \|f\|_\mathcal{H} + c'' t_n \|f\|_n \tag{50}$$

with probability at least $1 - c_1' \exp\left(-c_2' n \frac{t_n^2}{\sigma_v^2}\right)$.

**Proof.** The proof for bound (49) follows the proof for Lemma 13.2 in Wainwright (2015). To show (49), we let

$$\mathcal{A}(u) := \left\{ \exists f \in \bar{\mathcal{F}} \cap \{\|f\|_n \geq u\} : \frac{1}{n}\left|\sum_{i=1}^n v_i f(Z_i)\right| \geq 2 \|f\|_n u \right\}$$

where $u \geq r_n$. Suppose that there exists some $f \in \bar{\mathcal{F}}$ with $\|f\|_n \geq u$ such that

$$\frac{1}{n}\left|\sum_{i=1}^n v_i f(Z_i)\right| \geq 2 \|f\|_n u. \tag{51}$$

Let the function $\tilde{f} := \frac{u}{\|f\|_n} f$ and note that $\left\|\tilde{f}\right\|_n = u$. Since $f \in \bar{\mathcal{F}}$ and $u \leq \|f\|_n$ by construction, $\tilde{f} \in \bar{\mathcal{F}}$ under the star-shaped assumption. Consequently, whenever the event $\mathcal{A}(u)$ is true so that there exists a function $f$ satisfying inequality (51), then there exists a function $\tilde{f} \in \bar{\mathcal{F}}$ with $\|\tilde{f}\|_n = u$ such that

$$\frac{1}{n}\left|\sum_{i=1}^n v_i \tilde{f}(Z_i)\right| = \frac{u}{\|f\|_n}\frac{1}{n}\left|\sum_{i=1}^n v_i f(Z_i)\right| \geq 2u^2.$$

To summarize, we have established

$$\mathbb{P}\left[\mathcal{A}(u) \mid \{Z_i\}_{i=1}^n\right] \leq \mathbb{P}\left[A_n(u) \geq 2u^2 \mid \{Z_i\}_{i=1}^n\right],$$

where

$$A_n(u) := \sup_{\tilde{f} \in \Omega(u; \mathcal{F})} \frac{1}{n}\left|\sum_{i=1}^n v_i \tilde{f}(Z_i)\right|$$

with $\Omega(u; \mathcal{F}) = \left\{f \in \bar{\mathcal{F}} : \|f\|_n \leq u\right\}$. Conditioning on $\{Z_i\}_{i=1}^n$, note that under our assumptions, for each fixed $\tilde{f}$, the variable $\frac{1}{n}\sum_{i=1}^n v_i \tilde{f}(Z_i)$ is zero-mean sub-Gaussian. Then by standard tail bounds, we have

$$\mathbb{P}\left[A_n(u) \geq \mathbb{E}_v\left[A_n(u) \mid \{Z_i\}_{i=1}^n\right] + b \mid \{Z_i\}_{i=1}^n\right] \leq c_0 \exp\left(-c\frac{nb^2}{u^2\sigma_v^2}\right).$$

Consequently, for any $b = u^2$,

$$\mathbb{P}\left[A_n(u) \geq \mathbb{E}_v\left[A_n(u) \mid \{Z_i\}_{i=1}^n\right] + u^2 \mid \{Z_i\}_{i=1}^n\right] \leq c_0 \exp\left(-c\frac{nu^2}{\sigma_v^2}\right). \tag{52}$$

Finally, note that we have $\mathbb{E}_v\left[A_n(u) \mid \{Z_i\}_{i=1}^n\right] = \mathcal{G}_n(u; \mathcal{F})$. Since $u \geq r_n$ and the function $t \mapsto \frac{\mathcal{G}_n(t; \mathcal{F})}{t}$ is non-increasing (by Lemma A8), we obtain

$$\frac{\mathcal{G}_n(u; \mathcal{F})}{u} \leq \frac{\mathcal{G}_n(r_n; \mathcal{F})}{r_n} \leq \frac{r_n}{2},$$

where the last inequality uses the definition of $r_n$. Putting everything together, we have established that $\mathbb{E}_v\left[A_n(u) \mid \{Z_i\}_{i=1}^n\right] \leq \frac{ur_n}{2}$. Combined with the bound (52), we have

$$\mathbb{P}\left[A_n(u) \geq 2u^2 \mid \{Z_i\}_{i=1}^n\right] \leq \mathbb{P}\left[A_n(u) \geq ur_n + u^2 \mid \{Z_i\}_{i=1}^n\right]$$
$$\leq c_1 \exp\left(-c_2 n \frac{u^2}{\sigma_v^2}\right)$$

where the inequality uses the fact that $u^2 \geq ur_n$. This proves (49).

The second bound (50) follows from Lemma 1 in Raskutti et al. (2012). If $\|f\|_{\mathcal{H}} \lesssim 1$ and $\|f\|_n \lesssim t_n$, then $\left|\frac{1}{n}\sum_{i=1}^n v_i f(Z_i)\right| \leq c_3 t_n^2$ with probability at least $1 - c_1' \exp\left(-c_2' n \frac{t_n^2}{\sigma_v^2}\right)$. □

**Lemma A6** (LRE condition). Suppose Assumptions 1, 3, and (33) hold. Then, for any

$$\Delta_j \in \mathbb{C}(J(\pi_j)) := \left\{\Delta \in \mathbb{R}^{p-1} : \left|\Delta_{J(\pi_j)^c}\right|_1 \leq 3\left|\Delta_{J(\pi_j)}\right|_1\right\} \tag{53}$$

$[J(\pi_j)$ is the support of $\pi_j]$ and every $j = 1, ..., p$, we have

$$\Delta_j^T \frac{\hat{X}_{-j}^T \hat{X}_{-j}}{n} \Delta_j \geq \kappa_1 \|\Delta_j\|_2^2 \tag{54}$$

with probability at least $1 - c\exp\left(-c' n t_n^2 + c'' \log p\right)$ for any $t_n \geq r_n$, where $\kappa_1 = c_0' \Lambda_{\min}^2(\Sigma)$ for some universal constant $c_0' > 0$.

**Proof.** By elementary inequalities, we have

$$\left|\Delta_j^T \frac{\hat{X}_{-j}^T \hat{X}_{-j}}{n} \Delta_j\right| \geq \left|\Delta_j^T \frac{\tilde{X}_{-j}^T \tilde{X}_{-j}}{n} \Delta_j\right| - \left|\Delta_j^T \left(\frac{\tilde{X}_{-j}^T \tilde{X}_{-j} - \hat{X}_{-j}^T \hat{X}_{-j}}{n}\right) \Delta_j\right|$$

$$\geq \left|\Delta_j^T \frac{\tilde{X}_{-j}^T \tilde{X}_{-j}}{n} \Delta_j\right| - \left\|\frac{\tilde{X}_{-j}^T \tilde{X}_{-j} - \hat{X}_{-j}^T \hat{X}_{-j}}{n}\right\|_\infty \|\Delta_j\|_1^2$$

$$\geq \left|\Delta_j^T \frac{\tilde{X}_{-j}^T \tilde{X}_{-j}}{n} \Delta_j\right| - \left(\left\|\frac{\tilde{X}_{-j}^T(\hat{X}_{-j} - \tilde{X}_{-j})}{n}\right\|_\infty + \left\|\frac{(\hat{X}_{-j} - \tilde{X}_{-j})^T \hat{X}_{-j}}{n}\right\|_\infty\right) \|\Delta_j\|_1^2$$

$$\geq \left|\Delta_j^T \frac{\tilde{X}_{-j}^T \tilde{X}_{-j}}{n} \Delta_j\right| - \left\|\frac{\tilde{X}_{-j}^T(\hat{X}_{-j} - \tilde{X}_{-j})}{n}\right\|_\infty \|\Delta_j\|_1^2$$

$$- \left\|\frac{(\hat{X}_{-j} - \tilde{X}_{-j})^T \tilde{X}_{-j}}{n}\right\|_\infty \|\Delta_j\|_1^2 - \left\|\frac{(\hat{X}_{-j} - \tilde{X}_{-j})^T (\hat{X}_{-j} - \tilde{X}_{-j})}{n}\right\|_\infty \|\Delta_j\|_1^2.$$

We first bound $\left\|\frac{\tilde{X}_{-j}^T(\hat{X}_{-j} - \tilde{X}_{-j})}{n}\right\|_\infty$. Note that in terms of the $(l, l')$th element of the matrix $\frac{\tilde{X}_{-j}^T(\hat{X}_{-j} - \tilde{X}_{-j})}{n}$, for any $t_n \geq r_n$, we have

$$\left|\frac{1}{n}\tilde{X}_l^T(\hat{X}_{l'} - \tilde{X}_{l'})\right| = \left|\frac{1}{n}\sum_{i=1}^n \tilde{X}_{il}\left[\hat{f}_{l'}(Z_i) - f_{l'}(Z_i)\right]\right|.$$

Lemma A5 and (26) imply that, for any $t_n \geq r_n$,

$$\max_{l,l'}\left|\frac{1}{n}\tilde{X}_l^T(\hat{X}_{l'} - \tilde{X}_{l'})\right| \leq c_0 t_n^2 \tag{55}$$

with probability at least $1-c_1\exp\left(-c_2 nt_n^2+c_3\log p\right)$. To bound the term $\left\|\frac{(\hat{X}_{-j}-\tilde{X}_{-j})^T(\hat{X}_{-j}-\tilde{X}_{-j})}{n}\right\|_\infty$, we have

$$\left\|\frac{(\hat{X}_{-j}-\tilde{X}_{-j})^T(\hat{X}_{-j}-\tilde{X}_{-j})}{n}\right\|_\infty \leq \max_{l'}\left\|\hat{f}_{l'}(Z)-f_{l'}(Z)\right\|_n^2 \leq c_0' t_n^2 \qquad (56)$$

with probability at least $1-c_1'\exp\left(-c_2' nt_n^2+c_3'\log p\right)$.

To bound $\left|\Delta_j^T \frac{\tilde{X}_{-j}^T \tilde{X}_{-j}}{n}\Delta_j\right|$ from below, we apply the second result in Lemma A4; since this result holds for all $\Delta_j\in\mathbb{R}^{p-1}$, we can specialize it to any $\Delta_j$ in the restricted sets specified in (53). Putting everything above together and choosing $t_n=r_n\vee\sqrt{\frac{\log p}{n}}$ yields

$$\Delta_j^T \frac{\hat{X}_{-j}^T \hat{X}_{-j}}{n}\Delta_j \geq c_1\Lambda_{\min}^2(\Sigma)\|\Delta_j\|_2^2 - c_2\left(\frac{\log p}{n}\vee r_n^2\right)\|\Delta_j\|_1^2 \qquad (57)$$

with probability at least $1-c\exp\left(-c' nt_n^2+c''\log p\right)$. As a result, the cone condition $\|\hat{\pi}_j-\pi_j\|_1 \leq 4\sqrt{s_j}\|\hat{\pi}_j-\pi_j\|_2$ implied by (53) along with condition (33) yields (54). □

**Lemma A7** (Upper bound on $\max_j\left\|\frac{1}{n}\hat{X}_{-j}^T u_j\right\|_\infty$). Suppose $(\|\pi_j\|_1\vee 1)r_n^2=O\left(\sqrt{\frac{\log p}{n}}\right)$, Assumptions 1 and 3 hold. Then, we have $\max_{j=1,\ldots,p}\left\|\frac{1}{n}\hat{X}_{-j}^T u_j\right\|_\infty \precsim \max_{l,j}\mathbb{E}\left[\frac{1}{n}\tilde{X}_l^T\left(\tilde{X}_j-\tilde{X}_{-j}\pi_j\right)\right]+\sqrt{\frac{\log p}{n}}$ with probability at least $1-c_1'\exp\left(-c_2'\log p\right)-c_3'\exp\left(-c_4' nt_n^2+c_5'\log p\right)$.

**Proof.** Recall the definition of $u_j$ in (36) and

$$\hat{X}_{-j}=\tilde{X}_{-j}-\hat{\mathbb{E}}\left(X_{-j}|Z\right)+\mathbb{E}\left(X_{-j}|Z\right)$$

where $\tilde{X}_{-j}=X_{-j}-\mathbb{E}\left(X_{-j}|Z\right)$. For any $l\neq j$, after expanding $\left|\frac{1}{n}\hat{X}_l^T u_j\right|$, we note that in order to bound $\max_j\left\|\frac{1}{n}\hat{X}_{-j}^T u_j\right\|_\infty$, we need to bound $\max_{l'}\left\|\hat{f}_{l'}(Z)-f_{l'}(Z)\right\|_n^2$, $\max_{l,l'}\left|\frac{1}{n}\tilde{X}_l^T(\hat{X}_{l'}-\tilde{X}_{l'})\right|$ ($l\,l'=1,\ldots,p,\,l\neq j$), $\max_{l=1,\ldots,p,\,l\neq j}\left|\frac{1}{n}\eta_j^T(\hat{X}_l-\tilde{X}_l)\right|$, and $\max_{l=1,\ldots,p,\,l\neq j}\left|\frac{1}{n}\tilde{X}_l^T\eta_j\right|$. In particular, the first two terms are bounded in (56) and (55); the fourth term can be bounded by a standard sub-Exponential concentration inequality; for the third term, we evoke (46) and the argument that is used to bound (55), which yields

$$\max_{l\neq j}\left|\frac{1}{n}\eta_j^T(\hat{X}_l-\tilde{X}_l)\right|\leq c_0'' t_n^2$$

with probability at least $1-c_1''\exp\left(-c_2'' nt_n^2+c_3''\log p\right)$. Putting the pieces together and choosing $t_n=r_n\vee\sqrt{\frac{\log p}{n}}$, the claim in Lemma A7 follows. □

**Lemma A8.** Let the shifted function class $\bar{\mathcal{F}}$ be star-shaped. Then the function $t\mapsto\frac{\mathcal{G}_n(t;\mathcal{F})}{t}$ is non-increasing on the interval $(0,\infty)$; as a result, the *critical inequality* has a smallest positive solution (the *critical radius*).

**Remark.** This is Lemma 13.1 from Wainwright (2015).

**Lemma A9**: Let $N_n(t;\Omega(\tilde{r}_n;\mathcal{F}))$ denote the $t$−covering number of the set

$$\Omega(\tilde{r}_n;\mathcal{F})=\left\{f\in\bar{\mathcal{F}}:|f|_n\leq\tilde{r}_n\right\}$$

in the $\mathcal{L}^2(\mathbb{P}_n)$ norm. Then the smallest positive solution (the *critical radius*) to the *critical inequality* is bounded above by any $\tilde{r}_n \in (0, \sigma^\dagger]$ such that

$$\frac{c'}{\sqrt{n}} \int_{\frac{\tilde{r}_n^2}{4\sqrt{2}\sigma^\dagger}}^{\tilde{r}_n} \sqrt{\log N_n(t; \Omega(\tilde{r}_n; \mathcal{F}))} dt \leq \frac{\tilde{r}_n^2}{4\sigma^\dagger}.$$

**Remark**. This result is established by van der Vaart and Wellner (1996), van de Geer (2000), Barlett and Mendelson (2002), Koltchinski (2006), Wainwright (2015), etc.

**Lemma A10** (Approximately sparse $\pi_j$). Suppose Assumptions 1, 2(i) regarding $\Lambda_{\min}^2(\Sigma)$, 3, and 6 hold. If $(\|\pi_j\|_1 \vee 1) r_n^2 = O\left(\sqrt{\frac{\log p}{n}}\right)$, $\lambda_j \succsim \sqrt{\frac{\log p}{n}}$, and

$$\|\pi_j\|_1 \left(\frac{\log p}{n} \vee r_n^2\right) \leq c\sqrt{\frac{\log p}{n}} \Lambda_{\min}^2(\Sigma) \tag{58}$$

for some sufficiently small constant $c > 0$, then,

$$\max_j \|\hat{\pi}_j - \pi_j\|_2 = O_p\left(\lambda_j \sqrt{s_j} + \sqrt{\lambda_j \left\|\pi_{j,S_{\mathcal{T}_j}^c}\right\|_1}\right), \tag{59}$$

$$\max_j \|\hat{\pi}_j - \pi_j\|_1 = O_p\left(\lambda_j s_j + \sqrt{\lambda_j s_j \left\|\pi_{j,S_{\mathcal{T}_j}^c}\right\|_1} + \left\|\pi_{j,S_{\mathcal{T}_j}^c}\right\|_1\right). \tag{60}$$

**Proof**. Let $\hat{\Delta}_j = \hat{\pi}_j - \pi_j$. The basic inequality and the choice of $\lambda_j$ yield that $\left\|\hat{\Delta}_{j,S_{\mathcal{T}_j}^c}\right\|_1 \leq 3\left\|\hat{\Delta}_{j,S_{\mathcal{T}_j}}\right\|_1 + 4\left\|\pi_{j,S_{\mathcal{T}_j}^c}\right\|_1$. Consequently,

$$\left\|\hat{\Delta}_j\right\|_1 \leq 4\left\|\hat{\Delta}_{j,S_{\mathcal{T}_j}}\right\|_1 + 4\left\|\pi_{j,S_{\mathcal{T}_j}^c}\right\|_1 \leq 4\sqrt{s_j}\left\|\hat{\Delta}_j\right\|_2 + 4\left\|\pi_{j,S_{\mathcal{T}_j}^c}\right\|_1. \tag{61}$$

Moreover, we have

$$\sum_{l=1,\ldots,p, l \neq j} |\pi_{jl}| \geq \sum_{l \in S_{\mathcal{T}_j}} |\pi_{jl}| \geq \underline{\tau}_j s_j$$

and therefore $s_j \leq \underline{\tau}_j^{-1} \|\pi_j\|_1$. Putting the pieces together yields

$$\left\|\hat{\Delta}_j\right\|_1 \leq 4\sqrt{\underline{\tau}_j^{-1} \|\pi_j\|_1} \left\|\hat{\Delta}_j\right\|_2 + 4\left\|\pi_{j,S_{\mathcal{T}_j}^c}\right\|_1. \tag{62}$$

Therefore, for any vector $\hat{\Delta}_j$ satisfying (62), applying (57) yields

$$\hat{\Delta}_j^T \frac{\hat{X}_{-j}^T \hat{X}_{-j}}{n} \hat{\Delta}_j \geq \left\|\hat{\Delta}_j\right\|_2^2 \left\{c_1 \Lambda_{\min}^2(\Sigma) - c_2 \underline{\tau}_j^{-1} \|\pi_j\|_1 \left(\frac{\log p}{n} \vee r_n^2\right)\right\} - c_3 \left\|\pi_{j,S_{\mathcal{T}_j}^c}\right\|_1^2 \left(\frac{\log p}{n} \vee r_n^2\right). \tag{63}$$

With the choice of

$$\left(\Lambda_{\min}^2(\Sigma)\right)^{-\frac{1}{2}} \left\|\pi_{j,S_{\mathcal{T}_j}^c}\right\|_1 \sqrt{\frac{\log p}{n} \vee r_n^2} \asymp \delta^*,$$

under condition (58), we have

$$\left|\hat{\Delta}_j^T \frac{\hat{X}_{-j}^T \hat{X}_{-j}}{n} \hat{\Delta}_j\right| \geq c_4 \Lambda_{\min}^2(\Sigma) \left\|\hat{\Delta}_j\right\|_2^2$$

for any $\hat{\Delta}_j$ such that $\left\|\hat{\Delta}_j\right\|_2 \geq \delta^*$. We can then apply Theorem 1 in Negahban, et. al (2010) to obtain (59). By (61), we also obtain (60). $\square$

**Proof of Theorem 2**

For Theorem 2, we apply (59) and (60) when proving Lemma A3 and obtain

$$
\begin{aligned}
\frac{1}{\hat{\tau}_j^2} - \frac{1}{\tau_j^2} &= O_p\left(\max_j B_{2j}^2\right) + O_p\left(\max_j \left(\|\pi_j\|_1 \vee 1\right) \sqrt{\frac{\log p}{n}}\right) \\
&+ O_p\left(\max_j \left(\|\pi_j\|_1^2 \vee 1\right) \left(r_n^2 \vee \frac{\log p}{n}\right)\right), \\
\left\|\hat{\Theta}_j - \Theta_j\right\|_1 &= O_p\left(\max_j B_{1j}\right) + O_p\left(\max_j \left(\|\pi_j\|_1^2 \vee \|\pi_j\|_1\right) \sqrt{\frac{\log p}{n}}\right) \\
&+ O_p\left(\max_j \left(\|\pi_j\|_1^3 \vee \|\pi_j\|_1\right) \left(r_n^2 \vee \frac{\log p}{n}\right)\right), \\
\left\|\hat{\Theta}_j - \Theta_j\right\|_2 &= O_p\left(\max_j B_{2j}\right) + O_p\left(\max_j \left(\|\pi_j\|_1 \|\pi_j\|_2 \vee \|\pi_j\|_2\right) \sqrt{\frac{\log p}{n}}\right) \\
&+ O_p\left(\max_j \left(\|\pi_j\|_1^2 \|\pi_j\|_2 \vee \|\pi_j\|_2\right) \left(r_n^2 \vee \frac{\log p}{n}\right)\right), \\
\left|\hat{\Theta}_j \mathbb{E}\left(\tilde{X}_i^T \tilde{X}_i\right) \hat{\Theta}_j^T - \Theta_{j,j}\right| &= O_p\left(\max_j B_{2j}^2\right) + O_p\left(\max_j \left(\|\pi_j\|_1^2 \|\pi_j\|_2^2 \vee \|\pi_j\|_2^2\right) \frac{\log p}{n}\right) \\
&+ O_p\left(\max_j \left(\|\pi_j\|_1^4 \|\pi_j\|_2^2 \vee \|\pi_j\|_2^2\right) \left(r_n^2 \vee \frac{\log p}{n}\right)^2\right) \\
&+ O_p\left\{\max_j \left(\|\pi_j\|_1 \vee 1\right) \sqrt{\frac{\log p}{n}} + \max_j \left(\|\pi_j\|_1^2 \vee 1\right) \left(r_n^2 \vee \frac{\log p}{n}\right)\right\}, \\
\left\|\hat{\Theta}_j \frac{\hat{X}^T \hat{X}}{n} - e_j\right\|_\infty &= O_p\left(\max_j B_{1j}\right) + O_p\left(\max_j \left(\|\pi_j\|_1^2 \vee \|\pi_j\|_1\right) \sqrt{\frac{\log p}{n}}\right) \\
&+ O_p\left(\max_j \left(\|\pi_j\|_1^3 \vee \|\pi_j\|_1\right) \left(r_n^2 \vee \frac{\log p}{n}\right)\right).
\end{aligned}
$$

Now, we adopt the same argument as in the proof for Theorem 1 with the following minor differences: in showing the remainder terms $E_0 - E_4$ and $E_3' - E_6'$ are $o_p\left(\frac{1}{\sqrt{n}}\right)$ as well as that $\left|\hat{\Theta}_j \frac{\hat{X}^T \hat{X}}{n} \hat{\Theta}_j^T - \Theta_j \mathbb{E}\left(\tilde{X}_i^T \tilde{X}_i\right) \Theta_j^T\right| = o_p(1)$, we apply the above rates and Assumptions 4A-5A (which replace Assumptions 4-5). Putting these pieces together gives the claims in Theorem 2. $\square$

## 5 Simulations

In this section, we evaluate the performance of our methods with simulation experiments. To generate the full covariates $X$, we first generate $X_0$ from the $p$-dimensional normal distribution with mean 0 and variance $\Sigma_{X_0} = (\Sigma_{X_0,ij})_{i,j=1}^p$, which takes three different forms:

(S1) Independent: $\Sigma_{X_0} = I_p$;

(S2) AR(1): $\Sigma_{X_0,ij} = 0.5^{|i-j|}$;

(S3) Exchangeable/Compound Symmetric: $\Sigma_{X_0,ii} = 1$ and $\Sigma_{X_0,ij} = 0.5$ if $i \neq j$.

The covariates $\{Z_i\}_{i=1}^n$ are *i.i.d.* from $U[0, 2]$. To incorporate the dependence between $X$ and $Z$, we set $X_{i1} = X_{0,i1} + 3Z_i, X_{i2} = X_{0,i2} + 3Z_i^2, X_{i3} = X_{0,i3} - 3Z_i$ and $X_{ij} = X_{0,ij}, 1 \leq i \leq n, 4 \leq j \leq p$. The set of nonzero coefficients in $\beta_0$ is from a fixed realization of $s_0 = 3$ *i.i.d.* $U[0, 3]$. The active set is set to be $S_0 = \{1, 2, 3\}$. We consider two different non-linear functions $g_0$:

(G1) $g_0(z) = 1.5 \sin(2\pi z)$;

(G2) $g_0(z) = z^{10}(1-z)^4/B(11,5) + 4z^4(1-z)^{10}/B(5,11)$

where $B(\cdot, \cdot)$ denotes the beta distribution. The error terms are generated from a standard normal distribution. We estimate $\hat{X}$ by regressing $X$ on $Z$ with the smoothing spline procedure and $(\hat{\beta}, \hat{g})$ are obtained from (9). As a result, the debiased estimator in our simulations concerns (6). We do not test the performance of (5) with simulation experiments but expect it to behave similarly as (6). Similar to Zhang and Cheng (2017), the estimated variance $\hat{\sigma}_\varepsilon^2$ is calculated as follows:

$$\hat{\sigma}_\varepsilon^2 = \frac{\sum_{i=1}^n \left(Y_i - X_i\hat{\beta} - \hat{g}(Z_i)\right)^2}{n - \left\|\hat{\beta}\right\|_1}.$$

We set the tuning parameter $\mu = n^{-2/5}/10$ (Müller and van de Geer, 2015) and let $\lambda$ and $\lambda_j$ ($1 \leq j \leq p$) be calculated from the 10-fold cross validation (van de Geer, et al., 2014). Across all the simulations, we set the sample size $n = 100$ and the number of variables $p = 500$. Results in sections 5.1 and 5.2 are based on 100 replications, while those in section 5.3 are based on 500 replications.

## 5.1 Component-Wise Confidence Interval

Average coverage and average length of the intervals for individual coefficients corresponding to variables in either $S_0$ or $S_0^c$ are considered. Denote $\mathrm{CI}_j$ as a two-sided confidence interval for $\beta_j^0$. In Table 1, we report the empirical versions of

Avgcov $S_0 = s_0^{-1} \sum_{j \in S_0} \mathbb{P}(\beta_{0j} \in \mathrm{CI}_j)$;
Avgcov $S_0^c = (p - s_0)^{-1} \sum_{j \in S_0^c} \mathbb{P}(0 \in \mathrm{CI}_j)$;
Avglength $S_0 = s_0^{-1} \sum_{j \in S_0} \mathrm{length}(\mathrm{CI}_j)$;
Avglength $S_0^c = (p - s_0)^{-1} \sum_{j \in S_0^c} \mathrm{length}(\mathrm{CI}_j)$.

The results in Table 1 agree with our theoretical predictions. The average coverage probabilities of confidence intervals for $S_0^c$ are close to the nominal 95% level, while those for $S_0$ are slightly lower than 95%. The confidence intervals for $S_0^c$ are comparably narrower than those for $S_0$. We also notice that as the columns in $X_0$ become more correlated (so the inverse $\Theta$ of the Hessian becomes less sparse), the coverage performance becomes worse. This finding confirms our earlier comment (in Section 2) that the sparsity condition on the off diagonal elements of $\Theta$ plays a crucial role in the effectiveness of the debiased approach as this condition makes remainder terms like $\left(\hat{\Theta}_j - \Theta_j\right)\frac{1}{\sqrt{n}}\tilde{X}^T\varepsilon$ of order $o_p(1)$ in the asymptotic expansion of $\sqrt{n}\left(\tilde{b}_j - \beta_{0j}\right)$.

Table 1: Average coverage probabilities and lengths of confidence intervals at the 95% nominal level with 100 iterations; $n = 100, p = 500$

| Setup | Active set $S_0 = \{1,2,3\}$; Error $\varepsilon \sim N(0,1)$ | | | | | |
|---|---|---|---|---|---|---|
| Measure | S1, G1 | S2, G1 | S3, G1 | S1, G2 | S2, G2 | S3, G2 |
| Avgcov $S_0$ | 0.896 | 0.857 | 0.693 | 0.887 | 0.823 | 0.683 |
| Avglength $S_0$ | 0.802 | 0.812 | 0.807 | 0.798 | 0.789 | 0.827 |
| Avgcov $S_0^c$ | 0.953 | 0.955 | 0.963 | 0.953 | 0.955 | 0.963 |
| Avglength $S_0^c$ | 0.476 | 0.510 | 0.547 | 0.480 | 0.500 | 0.559 |

Table 2: Coverage probabilities and interval widths for the simultaneous confidence intervals based on the non-studentized (NST) and studentized (ST) test statistics with 100 iterations; $n = 100, p = 500$

| Setup | Active set $S_0 = \{1,2,3\}$; Error $\varepsilon \sim N(0,1)$ | | | | | |
|---|---|---|---|---|---|---|
| Measure | S1, G1 | S2, G1 | S3, G1 | S1, G2 | S2, G2 | S3, G2 |
| NST coverage | 0.95 | 0.74 | 0.72 | 0.96 | 0.86 | 0.74 |
| NST width | 1.05 | 1.09 | 1.12 | 1.06 | 1.06 | 1.12 |
| ST coverage | 0.88 | 0.94 | 0.87 | 0.82 | 0.91 | 0.83 |
| ST width | 0.88 | 0.96 | 1.04 | 0.87 | 0.93 | 1.04 |

## 5.2 Simultaneous Confidence Intervals

In Table 2, we present the coverage probabilities and interval widths for the simultaneous confidence intervals for $\beta_{0j}, 1 \leq j \leq p$. For each simulation run, we record whether the simultaneous confidence interval contains $\beta_{0j}$ for $1 \leq j \leq p$ and the corresponding interval width. Again, it is not surprising that the coverage probability is affected by the amount of correlations between the columns in $X_0$. Overall, both studentized and non-studentized method provide satisfactory coverage probability. When $\Sigma_{X_0}$ is the identity matrix, non-studentized method has better coverage; while when $\Sigma_{X_0}$ takes the form of S2 or S3, the performance of the studentized method is better.

## 5.3 Support Recovery

The major goal of this section is to identify signal locations of $\beta_0$ in a pre-specified set $G = \{1, 2, \ldots, p\}$, i.e. support recovery. Similarly as the procedure in Zhang and Cheng (2017), we take the signal set
$$\hat{\mathcal{S}}_0 = \{j \in \tilde{G} : |\tilde{b}_j| > \lambda_j^*\},$$

where $\lambda_j^* = \sqrt{2\hat{\omega}_{jj} \log(p)/n}$ and $\hat{\omega}_{jj} = \hat{\sigma}_\varepsilon^2 \hat{\Theta}_j \frac{\hat{X}^T \hat{X}}{n} \hat{\Theta}_j^T$. Note that similar arguments as Proposition 3.1 of Zhang and Cheng (2017) implies this support recovery procedure is consistent. To assess the performance, we consider the following similarity measure
$$d(\hat{\mathcal{S}}_0, \mathcal{S}_0) = \frac{|\hat{\mathcal{S}}_0 \cap \mathcal{S}_0|}{\sqrt{|\hat{\mathcal{S}}_0| \cdot |\mathcal{S}_0|}}.$$

Table 3 summarizes the mean and standard deviation of $d(\hat{\mathcal{S}}_0, \mathcal{S}_0)$ as well as the number of false positives (FP) and false negatives (FN) normalized by $\sqrt{|\hat{\mathcal{S}}_0| \cdot |\mathcal{S}_0|}$. When the amount of correlations between the columns in $X_0$ increases (as in S3), the false positive rates are comparably higher.

Table 3: The mean and standard deviation (SD) of $d(\hat{\mathcal{S}}_0, \mathcal{S}_0)$, and the numbers of false positives (FP) and false negatives (FN) with 500 iterations; $n = 100, p = 500$

| Setup | Active set $S_0 = \{1,2,3\}$; Error $\varepsilon \sim N(0,1)$ | | | | | |
|---|---|---|---|---|---|---|
| Measure | S1, G1 | S2, G1 | S3, G1 | S1, G2 | S2, G2 | S3, G2 |
| Mean | 0.96 | 0.97 | 0.94 | 0.97 | 0.97 | 0.94 |
| SD | 0.07 | 0.06 | 0.08 | 0.07 | 0.06 | 0.08 |
| FP | 0.09 | 0.06 | 0.13 | 0.09 | 0.07 | 0.12 |
| FN | 0.00 | 0.00 | 0.00 | 0.00 | 0.00 | 0.00 |